\documentclass[final]{siamltex}

\usepackage{algorithm,algorithmic}
\usepackage{amssymb}
\usepackage{MnSymbol} 

\usepackage{graphicx}
\usepackage[cmtip,all]{xy}
\usepackage{caption}
\newtheorem{context}[theorem]{Context}
\newtheorem{convention}[theorem]{Convention}

\newtheorem{remark}[theorem]{Remark}

\numberwithin{equation}{section}

\usepackage{color}
\usepackage{amsmath}
\usepackage{relsize} 
\usepackage{multicol}
\definecolor{MyDarkGreen}{rgb}{0,0.65,0}
\definecolor{MyDarkBlue}{rgb}{0,0,0.75}

\newcommand{\myqed}{\ensuremath{\square}}
\newcommand{\myqedhere}{\hfill \mbox{\raggedright \myqed}}
\newenvironment{proofsketch}{\par\emph{Proof (sketch).\ }}{\myqedhere}
\newenvironment{myproof}{\par\emph{Proof.\,}}{\myqedhere}
\newcommand{\overbar}[1]{\mkern 1.5mu\overline{\mkern-1.5mu#1\mkern-1.5mu}\mkern 1.5mu}
\def\inv{^{-1}}
\DeclareMathOperator{\mydeg}{deg} 
 
\DeclareMathOperator{\myspan}{span} 
\def\avgInt{\,\Xint-}

\def\XXint#1#2#3{{\setbox0=\hbox{$#1{#2#3}{\int}$}
     \vcenter{\hbox{$#2#3$}}\kern-.5\wd0}}
\def\Xint#1{\mathchoice
   {\XXint\displaystyle\textstyle{#1}}%
   {\XXint\textstyle\scriptstyle{#1}}%
   {\XXint\scriptstyle\scriptscriptstyle{#1}}%
   {\XXint\scriptscriptstyle\scriptscriptstyle{#1}}%
   \!\int}

\def\tplus{\!+\!}

\def\outcrowdingcap{boundary crowding cap\xspace}
\def\outcrowdingcaps{boundary crowding caps\xspace}
\def\outcrowding{boundary crowding\xspace}

\def\defining{\textbf}
\def\mention{\textit}
\def\capspd{\lambda}
\def\signal{\sigma}
\def\Signal{\mathcal{S}}
\def\Dt{\Delta t}
\def\Dtmax{\Dt_{\mathrm{max}}}
\def\Dtzero{\Dt_{\mathrm{zero}}}
\def\Dtstab{\Dt_{\mathrm{stable}}}
\def\Dtsafe{\Dt_{\mathrm{safe}}}
\def\Dtpos{\Dt_{\mathrm{pos}}}
\def\Dx{\Delta x}
\def\R{\mathcal{R}} 
\def\X{X}
\def\C{\mathcal{C}}
\def\CA{\overbar{\C}}
\def\B{\mathcal{B}}
\def\BA{\overbar{\B}}
\def\Bhat{\widehat\B} 
\def\BAhat{\widehat{\overbar{\B}}}
\def\sf{f}
\def\f{\mathbf{f}}
\def\h{h}
\def\vh{\underbar{h}}
\def\hA{\overbar{h}}

\def\U{U}
\def\vdU{\mathrm{d}\vU}
\def\hm{\h_{i-1/2}}
\def\hp{\h_{i+1/2}}
\def\Um{\U^-}
\def\Up{\U^+}
\def\Umph{\U^-_{i+1/2}}
\def\Ummh{\U^-_{i-1/2}}
\def\Upph{\U^+_{i+1/2}}
\def\Upmh{\U^+_{i-1/2}}
\def\Umm{\Ummh}
\def\Upp{\Upph}
\def\Ucm{\Upmh}
\def\Ucp{\Umph}
\def\Upm{\Ucm}
\def\Ump{\Ucp}

\def\Zcp{Z_+}
\def\Zcm{Z_-}
\def\dU{\mathrm{d}\U}
\def\du{\mathrm{d}\u}
\def\vdu{\mathrm{d}\vu}
\def\tU{\widetilde\U}
\def\Uavg{\overbar{U}}
\def\vUavg{\overbar{\vU}}
\def\opt{^\star}
\newcommand{\Uopt}[1][]{\ensuremath{{U\opt_{#1}}}}
\def\dual{^\ast} 
\def\PP{\mathbb{P}}
\def\M{M}
\def\MA{{\overbar{\M}}}
\def\MAinterval{\MAopt[\unitinterval]}
\def\unitinterval{{[0,1]}}
\newcommand{\WAopt}[1][]{\ensuremath{\WA{}\opt_{#1}}}
\newcommand{\MAopt}[1][]{\ensuremath{\MA{}\opt_{#1}}}

\newcommand{\Mopt}[1][]{\ensuremath{{M\opt_{#1}}}}
\newcommand\MAsph[2]{\ensuremath{{\MA{}_{\mysph{#1}}^{#2}}}}
\newcommand\MAbox[2]{\ensuremath{\MA{}_{\mybox{#1}}^{#2}}}
\newcommand\MAtri[2]{\ensuremath{{\MA{}_{\mytri{#1}}^{#2}}}}
\newcommand\MAstar[2]{\ensuremath{{\MA{}_{\ostar{#1}}^{#2,+}}}}
\def\MAinterval{\MA{}_{\unitinterval}}
\def\MAintervalk{\MAinterval^k}
\def\tight#1{\hbox{$#1$}}

\def\setof#1{\left\{ #1 \right\}}
\def\paren#1{\left( #1 \right)}
\def\bparen#1{\left[ #1 \right]}

\def\floor#1{\lfloor #1 \rfloor}
\def\tbracket#1{\![#1]\!}

\def\xm{{x^-}}
\def\xp{{x^+}}
\def\sp{s^+}
\def\sm{s^-}
\def\Sp{S^+}
\def\Sm{S^-}
\def\ssp{\dot{s}^+}
\def\ssm{\dot{s}^-}
\def\sSp{\dot{S}^+}
\def\sSm{\dot{S}^-}
\def\u{u}
\def\bu{\mathbf{u}}
\def\ubar{\overbar{u}}

\def\du{\mathrm{d}\u}

\def\vU{\underbar{U}}
\def\vu{\underbar{$u$}}
\def\uHLL{\u_{\hbox{\scriptsize HLL}}}

\def\energy{\mathcal{E}}
\def\mom{\mathbf{m}}
\def\calP{\mathcal{P}}
\def\thetabar{{\overbar{\theta}}}
\newcommand{\dampingcoef}[1][]{\ensuremath{\theta_{#1}}}

\def\Alpha{A}
\def\tAlpha{\widetilde\Alpha}
\def\tLambda{\widetilde\Lambda}
\def\thetaP{\dampingcoef[\calP]}
\def\thetaA{\dampingcoef[\Alpha]}

\def\and{\hbox{ and }}

\def\thmskip{\medskip}
\def\e{e}

\def\tn{\tilde n}

\def\x{\mathbf{x}}

\def\tu{\widetilde{u}}

\def\Umin{U_{\min}}
\def\Umax{U_{\max}}

\def\uhat{\widehat u}
\def\tu{\widetilde\u}
\def\rhoeps{\epsilon_\rho}
\def\peps{\epsilon_p}
\def\drho{\mathrm{d}\rho}

\def\nhat{\mathbf{\hat n}}

\def\vUm{\vU^-}
\def\vUQm{\vU^-_Q}
\def\vUQp{\vU^+_Q}
\def\vdUQm{\mathrm{d}\vU^-_Q}

\newcommand\ddt{\frac{\mathrm{d}}{\mathrm{d}t}}
\def\dt{\mathrm{d}_t}
\def\K{K}
\def\Q{Q}
\def\Y{Y}
\def\tK{{\widetilde K}}
\def\tX{{\widetilde X}}
\def\dK{{\partial K}}
\def\areaK{{|\partial K|}}
\def\volK{{|K|}}
\def\area{A}
\def\vol{V}
\def\dQ{\partial^Q}
\def\dQK{\dQ K}
\def\intK{\int_{\K}}
\def\QintK{\int^\mathrm{Q}_{\K}}
\def\intdK{\oint_{\dK}}
\def\QintdK{\oint^\mathrm{Q}_{\dK}}

\def\Qavg{\avgInt^\mathrm{Q}}

\def\avgK{\avgInt_{\K}}
\def\avgdK{\avgInt_{\dK}}
\def\QavgdK{\Qavg_{\dK}}
\def\dotp{{\,\boldsymbol\cdot\,}}
\def\Div{\nabla\dotp}
\def\dK{{\partial K}}
\def\next{^{n+1}}

\def\Uspace{\mathcal{V}}

\newcommand\Uplus[1][]{\ensuremath{\Uspace_{#1}^+}}
\newcommand\Uplushat[1][]{\ensuremath{\widehat{\Uspace}{}_{#1}^+}}

\def\d{\mathrm{d}}

\def\BAF{\BA^\mathrm{F}}

\def\WA{\overbar{\W}}

\newcommand{\Wopt}[1][]{\ensuremath{{W_{#1}}}}

\def\MX{\Mopt[\X]}

\def\W{W}

\def\maxspd{|\lambda|_{\max}}
\def\capspdA{\overbar{\capspd}}

\def\cs{\capspd}

\def\csA {\capspdA}
\def\csAF {\capspdA^\mathrm{F}}
\def\tR{\widetilde\R}
\def\RW{\R_{\W}}
\def\RM{\R^{\M}}
\def\RMA{\R^{\MA}}
\def\RA{\overbar \R} 
\def\RAWA{\RA{}_{\WA}}
\def\RAMA{\RA{}^{\MA}}

\def\RWA{\RAWA}
\def\RMA{\RAMA}
\def\Rhat{\widehat\R}

\def\reals{\mathbb{R}}

\def\realsgeZ{\reals_{\ge 0}}

\def\realsD{\reals^D}
\def\realsN{\reals^N}

\def\w{\mathrm{w}}
\def\wA{\overbar{\w}}
\def\wx{\w_{\x}}
\def\wAx{\wA_{\x}}
\def\twx{\widetilde\w_{\x}}
\def\what{\mathrm{\widehat w}}
\def\wbar{\mathrm{\overbar{w}}}
\def\wA{\wbar}
\def\nmax{n_{\max}}
\def\nmin{n_{\min}}
\usepackage{xspace}
\def\retentional{retentional\xspace}

\def\argmax{arg max\xspace}
\def\dr{\mathrm{d}r}
\def\Csum{\mathcal{C}} 
\def\Bsum{\mathcal{B}} 
\def\Cavg{\overbar \Csum} 
\def\Bavg{\overbar \Bsum} 

\def\PPD{\PP_D}
\def\PkD{\PP^k_D}

\def\dA{\mathrm{d}A}

\def\bigoh{\mathcal{O}}

\def\hAF{\hA^\mathrm{F}}
\def\csAF{\csA^\mathrm{F}}
\def\faces{\mathcal{F}}

\def\closure#1{\overbar{#1}}

\def\dP{\calP\dual}
\def\CP{\closure{\calP}}

\def\vbf{\underbar{$\mathbf{f}$}}

\def\Nk{N(k)}
\def\rhoavg{\overbar{\rho}}
\def\rhomin{\rho_\textrm{min}}

\def\uspd{|\bu|}

\def\ucap{u_\textrm{cap}}

\def\Ra{R_1}

\def\zsafety{\alpha_\mathrm{z}}

\def\uextra{u_\mathrm{extra}}
\def\Kbar{\overbar{\K}}

 \def\myleft({(}
 \def\myright){)}
\usepackage{comment}

\title{Outflow Positivity Limiting for Hyperbolic Conservation Laws.
  Part I: Framework and Recipe}

\author{Evan Alexander Johnson\thanks{KU Leuven, Department of Mathematics,
  Celestijnenlaan 200B box 2400, BE-3001 Heverlee, Belgium
  ({\tt e.alec.johnson@gmail.com}).}
        \and James A. Rossmanith\thanks{Iowa State University, Department of Mathematics,
396 Carver Hall, Ames, IA 50011, USA ({\tt rossmani@iastate.edu})}}

\begin{document}

\maketitle
\begin{abstract}
To support physically faithful simulation,
numerical methods for hyperbolic conservation laws
are needed that efficiently mimic the constraints satisfied by exact solutions,
including material conservation and positivity, while
also maintaining high-order accuracy and numerical stability.
Finite volume methods such as 
discontinuous Galerkin (DG) and weighted essentially non-oscillatory (WENO)
schemes allow efficient high-order accuracy while maintaining conservation.
Positivity limiters developed by Zhang and Shu and summarized
in [{\it Proc.\ R.\ Soc.\ A} {\bf 467}, 2752 (2011)] ensure a
minimum time step for which positivity of cell average quantities
is maintained without sacrificing conservation or formal
accuracy; this is achieved by linearly damping the deviation
from the cell average just enough to enforce a cell positivity
condition that requires positivity at boundary nodes and
strategically chosen interior points.

We assume that the set of positive states is convex; it follows that
positivity is equivalent to scalar positivity of a collection of affine
functionals.
Based on this observation, we generalize the technique of Zhang
and Shu to a framework that we call
outflow positivity limiting: 
First, enforce positivity at boundary nodes.
If wave speed desingularization is needed,
cap wave speeds at physically justified maxima
by using remapped states to calculate fluxes.
Second, apply linear damping again to cap the boundary average
of all positivity functionals
at the maximum possible (relative to the cell average)
for a scalar-valued representation positive in each mesh cell.
This be done by enforcing positivity of the \mention{retentional},
an affine combination of the cell average and the boundary average,
in the same way that Zhang and Shu would enforce
positivity at a single point (and with similar computational expense).
Third, limit the time step so that cell outflow
is less than the initial cell content.


We show that enforcing positivity at the interior points in Zhang
and Shu's method is actually a means of
capping boundary averages at the maximum possible for a positive solution.
Capping boundary averages allows computational
interior points to be chosen without sacrificing the guaranteed
positivity-preserving time step so as to optimize stabilization
benefits relative to computational expense, e.g.\ by choosing
points that coincide with nodal points of a DG scheme.

\end{abstract}




\section{Introduction}
\label{sec:introduction}
Consider a variable-coefficient hyperbolic scalar partial differential
equation (PDE) of the form
\begin{align}
  \label{conservationLaw}
  \partial_t\u(t,\x) +\Div\f(t,\x,\u) = 0
    \quad \hbox{ for } \quad \x\in\Omega, \, t\in\realsgeZ,
\end{align}
where $\Omega \subset \reals^D$ is the physical domain,
$t$ is the time coordinate,
$\Div$ represents the divergence with respect to $\x$,
$u:\realsgeZ\times \reals^D \rightarrow \reals$ is a conserved scalar quantity, and
the flux function $\f:\realsgeZ\times \reals^D\times \reals \rightarrow \reals^{D}$
is assumed differentiable.
Assume that solutions of $\eqref{conservationLaw}$ are positivity-preserving: that is,
if $u(0,\x) \ge 0$, then $u(t,\x) \ge 0$ for all $t>0$.
The context of this work will mostly be restricted to a single mesh cell
in the interior of the domain, so we do not concern ourselves with boundary conditions
for the PDE. 

A numerical method for equation \eqref{conservationLaw}
evolves a numerical solution $\U$ that approximates $\u$
and is represented by a finite number of degrees of freedom;
we say that the method is \emph{accurate} if the error
$\|\U-\u\|_\infty$ vanishes as the degrees of freedom are increased
in an appropriate way.

Accuracy is not the only goal of a numerical method, however.
One also seeks \emph{physicality} of numerical solutions.
If a physical condition is maintained
by the solution, then it is desirable for a
numerical method to mimic this by maintaining
a discrete version of this condition.
Such methods are referred to as \emph{mimetic} or \emph{conforming}.

For example, in a problem where material is conserved,
the solution $\u$ is restricted for all time
to a manifold of solutions that have the same total amount of material.
A numerical method that fails to satisfy a discrete
material conservation law can over time drift from this manifold,
resulting in unphysical behavior.  
Although one can correct for this by a global adjustment to the solution,
such an approach still can yield locally incorrect physics.
Specifically, solutions that fail to
satisfy a discrete \emph{local} conservation law result in simulated shocks
that travel at incorrect speeds (see e.g.\ pages 237--239 of LeVeque \cite{book:Le02}).

Similarly, in a problem where the solution should remain positive,
a numerical solution that fails to maintain a discrete
positivity condition could drift from the set of positive
solutions, resulting in an unphysical or even unstable solution.
Again, one could globally damp the deviation from the
global average.  But to ensure local physicality and stability,
positivity preservation should be enforced in a \emph{local} manner.

Finite volume methods are called \emph{conservative} because
they are designed to mimic the conservation
property of conservation laws.
Positivity-preserving methods are designed to
mimic the positivity-preserving property.
The challenge of positivity limiting is to 
design finite volume methods that are conservative,
positivity-preserving, high-order accurate, and numerically stable.

\begin{figure}[]
\fbox{\begin{minipage}{\linewidth}
\begin{center}
  \large{\textbf{Assumed properties and requirements}}
\end{center}
\vskip1ex
We assume a hyperbolic conservation law of the form
\begin{align}
  \label{hypsystemB}
  \partial_t\vu(t,\x) +\Div\vbf(t,\x,\vu) = 0,
\end{align}
where hyperbolicity means that
the flux Jacobian matrix
$(i,j)\mapsto \nhat\dotp\frac{\partial \f_i}{\partial u_j}$
is diagonalizable with real eigenvalues
and a full set of eigenvectors for any unit direction
vector $\nhat\in\reals^D$, where $\|\nhat\|=1$.
Each eigenvalue is a wave speed.
Valid physical solutions violate the differential form \eqref{hypsystemB}
at shock discontinuities but still satisfy
the integral form of conservation law \eqref{hypsystemB},
\begin{align}
  \label{integralConservationLaw}
  \mathrm{d}_t\intK\vu +\intdK\nhat\dotp\vbf = 0,
\end{align}
for arbitrary $\K$.
We assume that there is a convex set $\calP$
(satisfying $s\calP+(1-s)\calP\subset\calP\ \forall s\in[0,1]$)
of \mention{positive states} such that $\calP$ is an \mention{invariant domain}:
if the initial data $\vu(0,\x)$ is in $\calP$ then the physical solution
remains in $\calP$ for all time: $\vu(t,\x)\in\calP\ \forall t\ge 0$.
Furthermore, we assume that $\f\to 0$ on the boundary of $\calP$.
Positivity limiters yield a numerical method that is high-order accurate
for smooth solutions
and satisfies a discrete local conservation law like
\eqref{integralConservationLaw}
and a discrete local positivity condition $\intK\vu\in\calP$,
where $\K$ is restricted to a union of mesh cells.
\end{minipage}}
 \caption{Assumptions and requirements of positivity limiting
   (Section \ref{sec:introduction}).}
 \label{fundamentalAssumptions}
\end{figure}

Generalizing to hyperbolic systems, the fundamental assumptions and requirements
for positivity limiting are summarized in Figure \ref{fundamentalAssumptions}.

\subsection{Historical overview}

Consider a 1D conservation law
\begin{align*}
  \partial_t\u(t,x) + \partial_x f(\u) = 0
\end{align*}
which maintains the condition $\u\ge0$.

A finite volume method partitions the domain into intervals called
mesh cells. On a mesh cell of width $\Delta x$ centered at $x=x_i$
we denote the numerical cell average at time $t^n$ by $\Uavg_i^n$:
\begin{equation}
\label{eqn:1d_ave_update}
 \overbar{U}^{n}_i := \frac{1}{\Delta x} \int_{x_i - \frac{\Delta x}{2}}^{x_i + \frac{\Delta x}{2}} \, \U(t^n,x) \, dx.
\end{equation}
An Euler step for a method-of-lines finite volume method
(e.g.,\ DG or WENO) updates the cell average as follows:
\begin{gather}
\label{eqn:unlimited_rule}
 \Uavg\next_i = \Uavg^n_i - \tfrac{\Dt}{\Dx}\bparen{\h(\Umph,\Upph)-\h(\Ummh,\Upmh)}
\end{gather}
where $\Ummh$ and $\Upmh$ are approximate solution values on the left and right of cell interface $x_{i-1/2}$, respectively, 
and the numerical flux $h(\Um,\Up)$ is consistent with physical flux:
$h(u,u) = f(u)$.

In our abstract setting, what we know about $\f$ is that it is
defined so that physical solutions maintain positivity.
Therefore, to maintain positivity in the numerical scheme,
the numerical flux is defined in terms of the solution
to a physical problem.
This insight lead to the first general positivity-preserving scheme:
the \defining{Godunov scheme}\cite{God59}, which iterates the following:
\begin{enumerate}
\item In each mesh cell replace the solution with its average value.
\item Physically evolve the solution for a time step $\Dt$.
\end{enumerate}
If $\Dt$ is sufficiently short,
then the flux at each interface is given by the solution
to a \mention{Riemann problem}; this is guaranteed if
$\Dt\cs<\Dx$, where $\cs$ is the sum of left-going and right-going
signal speeds propagating from each interface in the physical solution;
see Figure \ref{defRiemannProblem}.
In this case, the Godunov scheme can be implemented by
equation \eqref{eqn:unlimited_rule} if the numerical flux
$h(\Um,\Up)$ is chosen to be the interface flux of the Riemann
problem with initial states $(\Um,\Up)$ and the solution
is assumed to be constant in each cell.
A simpler choice of $\h$ which maintains positivity
in update \eqref{eqn:unlimited_rule} is the HLL
numerical flux function \cite{article:HaLaLe83}. The HLL flux
is defined to account for the transfer of material implied
by the solution to the Riemann problem modified by averaging
in an interval that contains the signals emanating from the
interface. To show that HLL preserves positivity, Perthame and
Shu \cite{article:PerShu94} used a modified Godunov scheme: before
averaging in each mesh cell, the physically evolved solution is
first averaged in a set of nonoverlapping intervals each of which
contains the signals emanating from one of the interfaces.
The Godunov scheme is unfortunately only first-order accurate in space.


\begin{figure}[]
\fbox{\begin{minipage}{\linewidth}
\begin{center}
  \large{\textbf{Conservative local accuracy-preserving positivity limiting
    by linear damping of the deviation from the cell average}}
\end{center}
\vskip1ex
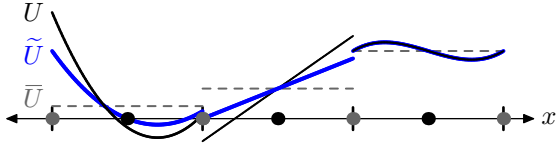
\begin{center}
\includegraphics{figures/riemann.6}
\end{center}
The positivity limiting framework introduced by Liu and Osher
\cite{article:LiuOsher96} achieves high-order accuracy while
maintaining positivity of the cell average $\Uavg$ in each mesh cell.
Prior to each time step the deviation from the cell average is linearly damped
until a cell positivity condition is satisfied; a 
time step is then calculated that preserves positivity of the cell average.
Zhang and Shu \cite{article:ZhShu10} reduce the cell positivity condition
from positivity at every point in the mesh cell
(i.e.\ for an infinite collection of point evaluation functionals) to positivity
at a \emph{finite} collection of \mention{positivity points}.
The limited solution is $\tU:=\Uavg+\theta(\U-\Uavg)$,
where $\theta\in[0,1]$ is just small enough so that $\tU$
is positive at each (black) interior point and at each
(gray) boundary node.
\end{minipage}}
\caption{Positivity limiting illustrated for 1D mesh cells with cubic polynomials.}
\label{depictionOfPositivityLimiting}
\end{figure}

To achieve high-order accuracy in space, the numerical flux function
$\h$ used in update \eqref{eqn:unlimited_rule}
needs to accurately approximate the physical flux $f$.
This will be the case if $\h$ is consistent with the physical
flux and if $\Um$ and $\Up$ are both high-order accurate
approximations to the exact solution at the interface.
Therefore, we assume in each cell a 
representation of the solution
of the form
\begin{equation}
      U^n_i(x) := \Uavg^{n}_i +  \dU^n_i,
\end{equation}
where $\Uavg^{n}_i$ is the cell average and
$\dU^n_i$ is a high-order correction
that is a polynomial of degree at most $k$.
The DG method evolves such a representation,
and WENO reconstructs such a representation with
each time step.  The flux function $h$ 
remains unchanged (i.e.\ is defined using the
Riemann problem with states $\Um$ and $\Up$.)

High-order accurate positivity limiters were developed
in the seminal works of Zhang \cite{zhang11thesis} and Zhang and Shu
\cite{article:ZhShu10}. An in-depth review of the history
of this method can be found in the excellent review article
\cite{article:ZhangShu11}.
As depicted in Figure \ref{depictionOfPositivityLimiting},
in each mesh cell
they linearly damp the high-order corrections just
enough to enforce positivity at a set of strategically
chosen positivity points.
For concreteness, we describe their method in detail
for the 1D scalar case.
Let $K=[x_i-\Dx/2,x_i+\Dx/2]$ denote a mesh cell,
and let
$\{s_j\}_{j=0}^{\tn}$
denote the points of
the Gauss-Lobatto quadrature rule that has just enough
points to exactly integrate polynomials of degree
at most $k$.


As depicted in figure \ref{depictionOfPositivityLimiting},
to limit $\U^n_i$, linearly damp $\dU^n_i$ just enough to
ensure that the solution is positive at each quadrature point.
In formulas:
\begin{equation}
\begin{split}
   \Umin := \min_{j} \, U^n_i(s_j), \quad
   \theta := \min \left\{ 1, \, \frac{\overbar{U}^{n}_i}{\overbar{U}^{n}_i - \Umin} \right\},
     \quad
  \tU^n_i(x) := \overbar{U}^{n}_i + \theta \, \text{d}U^{n}_i.
 \end{split}
   \label{eqn:limited_soln}
\end{equation}
We note that $\theta \in [0,1]$, and in particular, if $\Umin\ge 0$, then $\theta=1$.
Finally, the cell averages are updated by using the damped solution
in \eqref{eqn:unlimited_rule}:
\begin{equation}
\label{eqn:limited_update}
 \Uavg^{n+1}_i = \Uavg^{n}_i - \tfrac{\Delta t}{\Delta x} \, \left[ h\left(\tU^{-}_{i+1/2},\tU^{+}_{i+1/2}\right) - 
 h\left(\tU^{-}_{i-1/2},\tU^{+}_{i-1/2}\right) \right],
\end{equation}
where $\tU^{\pm}_{i-1/2}$ are the edge values computed from the
limited solution \eqref{eqn:limited_soln}.
We will see that this algorithm maintains accuracy
and guarantees the same positivity-preserving time step
that is guaranteed if positivity is enforced everywhere in the mesh cell.

Zhang and Shu have extended their positivity limiters 
to shallow water and gas dynamics
\cite{article:ZhangShu10rectmesh}, as well as
to higher-dimensional problems on both Cartesian and
unstructured grids \cite{article:ZhangShu10rectmesh,article:ZhangXiaShu12trimesh}.

As leveraged in \cite{article:ZhShu10}, a framework for
positivity limiting can also be used to implement a
maximum-principle-satisfying scheme. Suppose that the exact
(entropy) solution $\u(t,x)$ remains bounded in the interval
$[\Umin,\Umax]$.
By negation and translation of state space,
one can use scalar positivity limiting to enforce
that the affine functionals
$\Uavg^n_i-\Umin$ and $\Umax-\Uavg^n_i$ are both positive.

\subsection{Outflow positivity limiting framework}

\begin{figure}[]
\fbox{\begin{minipage}{\linewidth}
\begin{center}
{\large \textbf{Outflow capping}}
\end{center}
\begin{center}

\vskip-3ex
\begin{gather*}
 \begin{matrix}
    \begin{matrix}
      \mbox{\includegraphics[width=.47\linewidth]{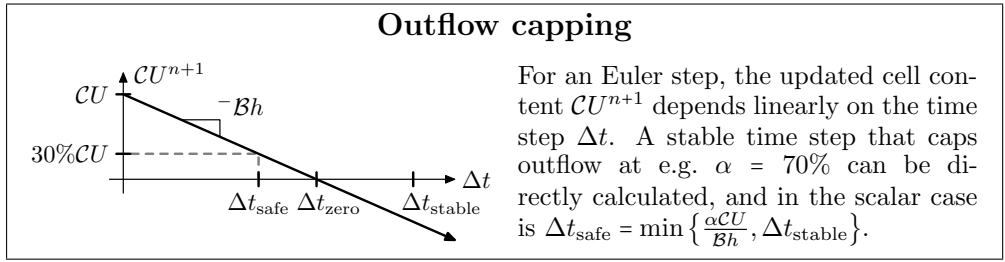}}
    \end{matrix}
    \quad
    \begin{matrix}
      \mbox{\parbox{.47\textwidth}{
       For an Euler step, the updated cell content $\C\U\next$
       depends linearly on the time step $\Dt$.
       A stable time step that caps outflow at e.g.\ $\alpha=70\%$
       can be directly calculated, and in the scalar
       case is $\Dtsafe=\min\setof{\frac{\alpha\C\U}{\B\h},\Dtstab}$.
       }}
    \end{matrix}
 \end{matrix}
\end{gather*}
\end{center}
\end{minipage}}
 \caption{Outflow capping in the scalar case.}
 \label{outflowCappingInTheScalarCase}
\end{figure}

The main contribution of this work is
an analysis and generalization of Zhang and Shu's positivity limiter
based on limiting the potential and actual outflow from each cell.
This yields a multi-stage framework we call {\it outflow positivity limiting}.
\begin{enumerate}

\item The first stage is \defining{point-wise positivity limiting},
which consists of enforcing positivity at
a set of \emph{positivity points} consisting of
\emph{boundary nodes} and \emph{interior points}.
Enforcing positivity at boundary nodes ensures that numerical
fluxes are defined with positive values. For shallow water and
gas dynamics, when positivity of the depth or density is enforced
at a boundary node, it is important to calculate fluxes with
remapped states in order to desingularize wave speeds.
Enforcing positivity at interior points may improve stability.

\item The second stage is \defining{boundary average limiting}. Like
Zhang and Shu, we maintain positivity of the cell average by
after each time step linearly damping the deviation from the
cell average just enough so that a cell positivity condition is
satisfied. The essential difference is that our cell positivity
condition directly caps the boundary average rather than requiring
positivity at points in the cell interior. Specifically, the
optimal version of our cell positivity condition requires the
boundary average to be no greater than the maximum boundary
average possible if the solution were everywhere positive and had
the same cell average.  Enforcing positivity at the optimal interior
points of Zhang and Shu is a means of boundary average limiting.

\item The third stage is \defining{outflow capping}. Rather than
restricting the time step to the guaranteed positivity-preserving
time step (and assuming that the estimated upper bound on wave
speeds is legitimate), we directly calculate a
maximal stable time step that caps cell outflow at e.g.\ 70\%.
The updated cell average varies linearly with the length
of an Euler time step, and therefore outflow capping can be
performed with the same limiting procedure that Zhang and
Shu use to enforce positivity at a positivity point.
See Figure \ref{outflowCappingInTheScalarCase}.
\end{enumerate}

\subsection{Benefits of the outflow positivity limiting framework}

This framework offers
simplicity, insight, flexibility, and computational efficiency.
Our framework is implied by a set of natural requirements,
e.g.\ to maintain positivity of the cell average
and to guarantee a stable, positivity-preserving time step.
We thereby decouple requirements and arrive
at a framework that is intrinsic to the requirements
and arguably optimal.
Outflow capping ensures that
knowledge of the fastest wave speed is not needed
to compute the positivity-preserving time step.
A minimum positivity-preserving time step is guaranteed
by positivity of the retentional, which is defined in terms of
the data that actually determines the evolution of the cell average
(i.e.\ boundary node data and the initial cell average),
and has a physical meaning (cell content minus maximum possible loss).
The proof that positivity is maintained is thus intuitively
obvious and works for spatially varying flux functions and
for arbitrary cell geometry and representation space.

Boundary average limiting can be used to guarantee a
positivity-preserving time step even if the data used at
the boundary nodes is not interpolated by a polynomial.
This is relevant in the gas-dynamics case, where it
becomes necessary to modify states at boundary nodes
before evaluating fluxes, e.g.\ to desingularize fluid velocity
when enforcing positivity of the density or to enforce
realizability or hyperbolicity in higher-moment gas-dynamic models.

\begin{figure}[]
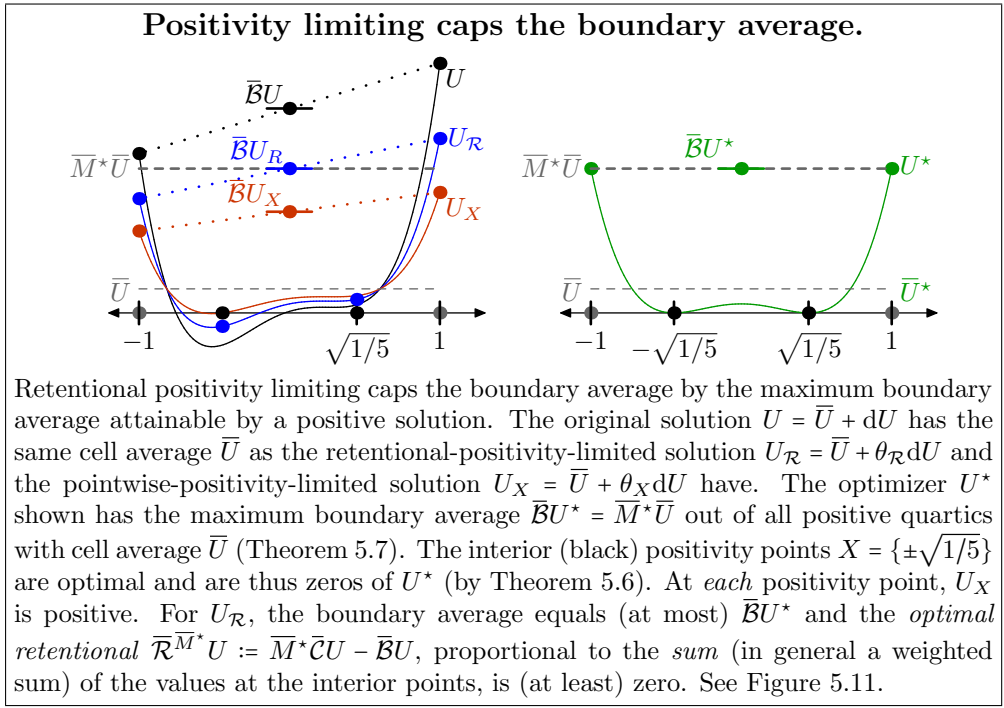

\fbox{\begin{minipage}{\linewidth}
\begin{center}
  \large{\textbf{Positivity limiting caps the boundary average.}}
\end{center}
\vskip1ex
\begin{center}
\includegraphics{figures/riemann.3}
\quad
\includegraphics{figures/riemann.4}
\end{center}
Retentional positivity limiting caps the boundary average by
the maximum boundary average attainable by a positive solution.
The original solution $\U=\Uavg+\dU$ has the same cell average $\Uavg$
as the retentional-positivity-limited solution $\U_\R=\Uavg+\theta_\R\dU$ and
the pointwise-positivity-limited solution $\U_\X=\Uavg+\theta_\X\dU$ have.
The optimizer $\Uopt$ shown has the maximum boundary average
$\BA\Uopt = \MAopt\Uavg$ out of all positive quartics with
cell average $\Uavg$ (Theorem \ref{optimizerExists}).
The interior (black) positivity points $\X=\{\pm\sqrt{1/5}\}$
are optimal and are thus zeros of $\Uopt$
(by Theorem \ref{optimalBoundaryWeightedQuadratureExists}).
At \emph{each} positivity point, 
$\U_\X$ is positive.
For $\U_\R$, the boundary average equals (at most) $\BA\Uopt$
and the \emph{optimal retentional}
$\RA^{\MAopt}\U:=\MAopt\CA\U-\BA\U$,
proportional to the \emph{sum} (in general a weighted sum)
of the values at the interior points, is (at least) zero.
See Figure \ref{mutualConfirmation}.
\end{minipage}}
\caption{Retentional positivity limiting for a one-dimensional mesh cell with a
quartic polynomial basis ($\MAopt=6$), compared with pointwise positivity limiting.}
\label{pointwiseVersusRetentionalPositivityLimiting}
\end{figure}

Capping the boundary average is enforced via positivity of a
single linear functional, the retentional, and is therefore no more
expensive than enforcing positivity at a single point.
%
For stabilization purposes, one may additionally enforce
positivity at interior points, but the choice of these points
can be optimized for requirements of efficiency, simplicity, or
stability; if positivity points are instead
optimized for maximizing the guaranteed
positivity-preserving time step, then enforcing
positivity at these points enforces the optimal cap on the
boundary average (as illustrated in Figure
\ref{pointwiseVersusRetentionalPositivityLimiting}.)
While we show that such optimal positivity
points always exist (see Figure \ref{mutualConfirmationFig}),
precisely isolating them can be non-trivial
(in general requiring the solution of linear programming
problems on a convex domain), and in the case of nodal DG it is
inefficient to check positivity at these non-nodal points if
they are relatively numerous. In contrast, outflow positivity
limiting merely requires an estimate of the maximum \outcrowding
$\MAopt$ realizable by a positive solution.  Note 
that precisely isolating $\MAopt$ also requires solving
linear programming problems on a convex domain.

\subsection{Parts}

This work is conceived of as the first in a three-part series.
This first part is analytical rather than computational
and avoids making claims dependent on
computational experiment; rather, we lay out a general
framework that incorporates previous work and that will guide
subsesquent computational investigation.
The second part will give a detailed justification of the
analytical results. The third part is planned to consider issues
that require computational simulation. In particular, the
stability consequences of enforcing positivity at Zhang and Shu's
points, at all nodal points, or only at boundary nodes, will be a
natural follow-up study, and results will likely be dependent on
the details of the choice of oscillation-suppressing limiters.
Since our framework generalizes that of Zhang and Shu, it has
in a sense already been computationally tested
by their work 
\cite{WangZhangShuNing12,zhang11thesis,article:ZhShu10,article:ZhangShu10rectmesh,article:ZhangShu11,article:ZhangShu10source,article:ZhangXiaShu12trimesh}.

\subsection{Section summary}

The outflow positivity limiting
framework is developed in subsequent sections of this paper.
In Section \ref{sec:dg_schemes} we describe how DG (or WENO)
updates the cell average.
In Section \ref{sec:framework} we show in the context of scalar
conservation laws that enforcing positivity
of the retentional guarantees that the cell average remains
positive for a time step inversely proportional to the 
\outcrowdingcap $\MA$ chosen in the definition of the retentional.
In Section \ref{sec:systems} we show that the same statement
holds for hyperbolic systems of conservation laws.
In Section \ref{sec:accuracy} we develop a theory of admissible
(accuracy-preserving) and optimal cell positivity conditions that
shows that, as long as the \outcrowdingcap is admissible
(i.e.\ no less than an optimal $\MAopt$), then
enforcing positivity of the retentional by linearly damping the deviation from the
cell average preserves the order of accuracy and can be accomplished by enforcing
positivity at strategically chosen interior points.
In Section \ref{CalculationOfWeights} we apply this theory of
cell positivity conditions to compute optimal or near-optimal
\outcrowdingcaps and optimal
interior points for typical cell geometries and polynomial orders.
In Section \ref{practicalApplication} we discuss practical
implementation issues regarding efficiency, robustness, and generalization,
with specific reference to shallow water and gas dynamics.
These issues include efficient and correct positivity checks, 
wave speed desingularization, stability benefits of interior positivity points,
stable optimal time stepping, multistage and local time stepping, 
and application to isoparametric mesh cells.
In Section \ref{Conclusion} we summarize the key results of this
work.  Boxed figures present a largely self-contained summary of the work.


\section{High-order discontinuous Galerkin schemes}
\label{sec:dg_schemes}
For simplicity, we discuss positivity-preserving limiters with
reference to the discontinuous Galerkin (DG) method. This does
not entail loss of generality, since a WENO finite volume method
can be thought of as a DG method that reconstructs the high-order
component of the representation prior to each time step. Standard
DG schemes discretize the weak form of hyperbolic conservation
law \eqref{conservationLaw}:
\begin{gather}
  \label{weakEvolution}
  \ddt\int_{\K} \u \, \varphi 
   = \int_{\K} \f\dotp \nabla \varphi - \oint_{\partial \K} \nhat\dotp\f \, \varphi,
\end{gather}
where $\K$ is an open subset of $\reals^D$, $\partial \K$ is the boundary of $\K$, $\nhat$ is an outward pointing unit
vector that at every point on $\partial \K$ is perpendicular to the boundary $\partial \K$, and $\varphi\in C^1_0$ is a test function.

Assume that $\Omega$ is the union of a finite set $\Omega_h$
of non-overlapping mesh cells.
We assume that the solution, when restricted to any mesh cell $\K\in\Omega_h$,
belongs to a polynomial space $\Uspace=\Uspace(\K)$ and is otherwise unconstrained
(e.g.\ by the requirement of continuity across mesh cell boundaries).
Often $\Uspace$ is $\PkD$,
the set of polynomials in $D$ variables of degree at most $k$.
More generally, we define $\mydeg(\Uspace)\ge 0$ to be the largest $k$
such that $\Uspace$ contains $\PkD$.
The discontinuous Galerkin method approximates the exact solution to \eqref{weakEvolution}
with a piecewise polynomial function:
\begin{equation}
\U\bigl|_\K=\sum_{i=0}^{\Nk} \U^{i} \, \varphi_{i} \, \in\Uspace,
\end{equation}
where, on each cell $K$, $\{\varphi^{j}\}_{j=0}^{\Nk}$ is a basis
for $\Uspace(\K)$ and $\{\varphi_{j}\}_{j=0}^{\Nk}$ is its co-basis.
The basis and co-basis are mutually orthonormal in the inner product
$ 
\frac{1}{\volK} \intK \varphi^{i} \varphi_{j} = \delta^i_j,
$ 
where $\volK$ is the measure of the mesh cell $K$.

In \defining{nodal DG}, the co-basis consists of a set of
point-evaluation functionals that evaluate the solution at the nodes,
and the solution is thus represented by nodal values.
In \defining{modal DG}, the basis is often chosen to be a sequence of
orthonormal polynomials of increasing degree.
As a consequence, the basis and co-basis are identical
and the coefficients of higher-order basis functions decay:
specifically, for a smooth function
$\u$, the $L^2$-projection coefficients $\intK \u \, \varphi^j$ decay
like $\bigoh(\Dx^{k'+1})$, where $\Dx$ is the
mesh cell diameter and $k'$ is largest such that 
$\PPD^{k'}\subset \myspan\setof{\varphi^i: i<j}$.

\begin{figure}[]
\fbox{\begin{minipage}{\linewidth}
The \defining{Riemann problem with
states $(\Um,\Up)$ and flux function $f$}
is defined to be the problem
\def\PDE{\mathrm{(PDE)}}
\begin{gather*}
  \PDE\quad
  \partial_t u(t,x) + \partial_x f(u)=0, \qquad
 \begin{gathered}
  \u(0,x) = 
  \left\{
    \begin{array}{l l}
      \Um & \hbox{ if } x<0,
   \\ \Up & \hbox{ if } x>0.
    \end{array}
  \right.
 \end{gathered}
\end{gather*}
Riemann problems (and thus their solutions) are by requirement
\emph{self-similar}, i.e.\ invariant
under dilation of space-time: $u(at,ax)=u(t,x)$ for any $a>0$ and $t\ge 0$.
So we can write $u(t,x) =: \uhat(x/t)$ for $t>0$. 
In hyperbolic problems information propagates at finite speed.
Define the \defining{signals}
$\sm(t)\le0\le\sp(t)$
emanating from the \mention{interface} $x=0$
so that $\signal:=[\sm,\sp]$
is the smallest interval such that
$u(t,\cdot)$ agrees with $u(0,\cdot)$ outside $\signal$.
By self-similarity, the \mention{signal speeds}
$\dot{s}^\pm$ are independent of $t$ and define
the \defining{right-going signal speed} $\ssp(\Um,\Up)$
and \defining{left-going signal speed} $\ssm(\Um,\Up)$.
The \mention{interface flux}
$f_0$ is defined to be the flux
at $x=0$ needed to account for the amount of material that
$u$ accumulates on each side of the interface; it is
independent of $t>0$ and equals 
$f(\uhat(0))$ unless $\tu$ happens to be discontinuous at $0$;
then subtracting
$\int_\xm^0\PDE$ from 
$\int_0^\xp\PDE$, applying the fundamental theorem of calculus,
and solving 
reveals that
$ 
  2 f_0 = f(\Um)+f(\Up)+
    \dt\big[\int_0^{\xp}\u(t,\cdot)-\int_{\xm}^0\u(t,\cdot)\big]
$ 
if $\xm<t\sm$ and $\xp>t\sp$.
\end{minipage}}
 \caption{Definition of the Riemann problem used to define numerical flux}
 \label{defRiemannProblem}
\end{figure}

The DG method solves a discretization of equation
\eqref{weakEvolution} in each mesh cell $\K$:
\def\QintK{\int^\mathrm{Q}_{\K}}
\begin{gather}
  \label{DGintEvolution}
  \ddt \intK\U\varphi^{j}  = \QintK \f(t,\x,\U) \dotp\nabla\varphi^{j} - \QintdK\h(U^-, U^+, \nhat) \, \varphi^{j},
\end{gather}
where $\nhat$ is the outward pointing unit normal to $\partial K$ and 
where $\h(\Um,\Up,\nhat)$ denotes a numerical flux function
where $\Um$ represents the state at the boundary node approached
from inside the cell and $\Up$ represents the state at
the same boundary node approached from outside the cell.
See Figure \ref{methodOfLinesFV}.
%
In equation \eqref{DGintEvolution} we have replaced exact integration by
quadrature rules; in particular,  $\QintK$ denotes a quadrature rule exactly equal
to $\intK$ for polynomials of degree at most $2k-1$ and
$\QintdK$ denotes a quadrature rule exactly equal to $\intdK$
for polynomials of degree at most $2k$.
If $\K$ is a polytope 
then the boundary quadrature rule can be
written as the sum over all faces of a Gaussian quadrature rule on each face.

\begin{figure}[]
\fbox{\begin{minipage}{\linewidth}
\begin{center}
  \large{\textbf{Framework of a method-of-lines finite volume Euler step}}
\end{center}
\vskip1ex
Integrating the hyperbolic differential equation
\begin{align*}
  \partial_t\u(t,\x) +\Div\f(t,\x,\u) = 0
\end{align*}
over a mesh cell $\K$ gives the integral form
\begin{align*}
  \ddt\intK\u +\intdK\nhat\dotp\f = 0.
\end{align*}
Making the replacements $\u\rightarrow\U\in\Uspace$,
$\intdK\rightarrow\B$, and
$\nhat\dotp\f\rightarrow\h$,
where
$\Uspace$ is a finite-dimensional representation space,
$\h$ is a numerical flux, and
$\B$ is a numerical boundary quadrature,
gives the ordinary differential equation of
a method-of-lines finite volume approximation
$ 
  \d_t\C\U = -\B\h,
$ 
where $\C:=\intK$.
An Euler step for this ODE is
\begin{gather*}
  \C\U\next = 
  \C\U - \Dt\B \h.
\end{gather*}
Typically $\Uspace$ when restricted to a mesh cell
is a polynomial representation space
containing all polynomials of degree at most $k$,
$\QintdK\h=\sum_{\x\in\Q}\omega_{\x}\h(\x)$
is a numerical quadrature rule with points $\Q\subset\dK$
and weights $\omega_{\x}>0$, and $\h(\x)$ is a numerical
interface flux for the Riemann problem with frozen flux function
$f:=\u\mapsto\nhat\dotp\f(t,\x,\u)$
and states $(\Um,\Up)$ (see Figure \ref{defRiemannProblem});
here $\Um(\x)$ is the value of $\U$ at $\x$ approached from
within $\K$ and $\Up(\x)$ is the value of $\U$ at $\x$ approached
from outside $\K$.
DG and WENO are two ways of updating the deviation of $\U$
from the cell average.
\end{minipage}}
 \caption{The context of this work as outlined in Section \ref{sec:dg_schemes}.}
 \label{methodOfLinesFV}
\end{figure}

An important property of the DG scheme is that it conserves the
total amount of $U$ from time step to time step; indeed, taking
$\varphi=1\in\Uspace$, equation
\eqref{DGintEvolution} becomes the conservation law
\begin{align}
  \label{FVintEvolution}
  \ddt\intK\U = -\QintdK\h.
\end{align}
%

The entire context of this work is concerned with maintaining
positivity in a single arbitrary mesh cell.
We therefore adopt a simplified notation that omits
explicit reference to the mesh cell $\K$.
We define the \defining{cell content}
to be the cell integral 
$\C := \U\mapsto\intK\U$
or its scale-invariant version, the cell average
$\CA := \U\mapsto\avgK\U$.
We define the \defining{boundary sum}
to be the 
boundary integral quadrature
$\B := \U\mapsto\QintdK\Um$
or its scale-invariant version, the boundary average quadrature
$\BA := \U\mapsto\QavgdK\Um$.
We use $\vol:=\volK$ to denote the ``volume'' of the mesh cell and
$\area:=\areaK$ to denote the ``area'' of its boundary.
With these conventions, the ODE \eqref{FVintEvolution}
is expressed as
\begin{align*}
  \d_t\C\U &= -\B\h,\qquad\hbox{i.e.,} \\
  \d_t\CA\U &= -\BA\hA,
\end{align*}
where $\hA:=\frac{\area}{\vol}\h$
is a scale-invariant version of the numerical flux,
and an explicit Euler step is expressed as
\begin{align}
  \label{explicitEuler}
  \C\U\next &= 
    \C\U - \Dt\B \h, \qquad \hbox{i.e.,}
  \\
  \label{explicitEulerA}
  \CA\U\next &= 
    \CA\U - \Dt\BA \hA.
\end{align}
This method-of-lines finite volume framework 
is summarized in Figure \ref{methodOfLinesFV}.

\section{Framework for outflow positivity limiting}
\label{sec:framework}

In this work \defining{positive} means \emph{nonnegative}
unless qualified with the adjective \emph{strictly}.
A \mention{positive combination} means a linear combination
with positive coefficients, and
a \mention{positive functional} is defined to be a
functional that is positive on positive functions.

A simple algorithm that maintains positivity of the cell average is
to repeat the following sequence:
\begin{enumerate}
 \item Assume that the cell average is positive.
 \item If necessary, linearly damp (rescale) the deviation
   from the cell average just enough so that a cell positivity
   condition is satisfied.  (See Figure \ref{depictionOfPositivityLimiting}.)
 \item Execute an Euler step for a stable time step that is just short
   enough so as to guarantee that positivity
   of the cell average is maintained (given that the
   cell positivity condition is satisfied).
\end{enumerate}
This algorithm was first introduced by Liu and Osher \cite{article:LiuOsher96}
and further developed by Zhang and Shu \cite{article:ZhShu10}.


Assuming that numerical fluxes are calculated by evaluating the flux
at the boundary node states (rather than at remapped states),
the cell positivity condition should at least require the solution to be
positive at the boundary nodes so that the Riemann problem used
to define the numerical flux function is well-defined.
To ensure that positivity limiting does not compromise order of accuracy,
it will be enough to require that the cell positivity condition
be satisfied if the initial data is already everywhere positive
(as we will verify in Theorem \ref{OptimalWeightPreservesAccuracy}).

Given the constraints of this framework, we ask:
\emph{What cell positivity condition and time step will
guarantee that $\avgK\U\next\ge 0$?}
If we are not allowed to modify a positive solution, then
the maximum time step for which we can hope to guarantee
positivity is the maximum such time step if the initial
data is positive.  We therefore ask:
  \emph{What is the maximum time step for which it
    is guaranteed that $\avgK\U\next\ge0$
    if $\U\ge 0$ in $\K$?}
We can assume that the outflow $\QintdK \h$ is
strictly positive, since otherwise the cell average
will remain positive for any positive time step.
In this case, setting $\C\U\next\ge0$ in equation
\eqref{explicitEuler} or \eqref{explicitEulerA} gives
\begin{gather}
  \label{ExactPositivityCondition}
  \Dt\inv \ge
     \Dtzero\inv
    := \frac{\B \h}{\C\U}
     = \frac{\BA \hA}{\CA\U}
\end{gather}
That is, the time $\Dtzero$ until strict positivity of the cell average is
violated is the ratio of the total integral of $U$ in the cell to the
net rate of flux out of the boundary. Clearly, to
guarantee positivity we need a constraint on the flux out of the
boundary. A constraint that extends naturally to systems is to
specify a cap on wave speeds.  Abstractly we impose that
\begin{gather}
  \label{hConstraint}
  \h \le \cs \Um,
\end{gather}
where we call $\cs$ the \defining{speed cap}; equivalently,
$ 
  \hA \le \csA \Um,
$ 
where $\csA:=\frac{\area}{\vol}\cs$
is a scale-invariant version of the speed cap.
This condition clearly holds for a scalar problem with a convex
flux function (e.g.\ Burgers' equation) if $\cs$ is the wave
speed, and in Theorem
\ref{speedCapTheorem} we show that for general systems $\cs$ can
be defined as a cap on the sum of incoming and outgoing signal speeds
at each boundary node.

\begin{remark}[\emph{omitting bars}]
\label{droppingBars}
\emph{
The form of the scale-invariant equations is identical
to the form of the non-invariant equations, so we can choose
whether to omit or retain bars when discussing general properties.
}
\end{remark}

\subsection{Boundary crowding caps and the retentional}
\label{capsAndRetentional}

Equations
\eqref{ExactPositivityCondition}
and
\eqref{hConstraint}
imply that
positivity is maintained if
\begin{gather}
  \label{positivityCondition}
  \Dt\inv \ge
      \cs\Bhat\U 
     = \csA\BAhat\U, 
\end{gather}
where we call
$\Bhat\U:= \frac{\B \U}{\C\U}$
or its scale-invariant version
$\BAhat\U:= \frac{\BA \U}{\CA\U}$
the \defining{\outcrowding.}
Therefore, we can guarantee a minimum positivity-preserving
time step by enforcing bounds on $\cs$ and $\Bhat\U$.
To maintain order of accuracy, these bounds need to be
physically justified.
Enforcing a cap on $\cs$ is briefly considered in
Section \ref{sec:desing}.  The focus of this paper
is on enforcing a justified cap on $\Bhat\U$.

Given that $\MA\ge\BAhat\U$,
equation \eqref{positivityCondition}
says that positivity is maintained if
\begin{gather}
  \label{weightFunctionalCondB}
  \Dt\inv \ge \Dtpos\inv := \csA\MA.
\end{gather}
To determine a justified $\MA$,
we use that physical solutions satisfy positivity:

\begin{definition}[$\MAopt$]
We define the \defining{positive solutions} to be
the set 
$\Uplus[\K]:=\setof{\U\in\Uspace: \U\ge0 \hbox{ in } \K \and \C\U>0}$
of representations positive and somewhere
nonzero in the mesh cell,
and we define the \defining{optimal \outcrowding cap}
(or \defining{optimal interior weight}) $\MAopt$ to be the
maximum \outcrowding over all positive solutions:
\begin{align}
  \label{MAoptDef}
  \MAopt[\K](\Uspace) &:= \sup_{\U\in\Uplus[\K]} \BAhat(\U).
\end{align}
\end{definition}
We can enforce that $\MA\ge\BAhat\U$ in exactly the same way
that positivity is enforced at positivity points if we define
this condition in terms of positivity of a linear functional:
\begin{definition}[retentional]
Enforcing a cap $\MA$ on $\BAhat\U=\frac{\BA\U}{\CA\U}$ is equivalent to 
enforcing positivity of the \defining{retentional}\,\footnote{
An abbreviated form of \emph{retention functional},
in the admittedly hokey mathematical tradition of using
``adjectivals'' as substantives and thus nouns.}
$
  \RAMA\U := \MA\CA\U-\BA\U,
$
or one of its rescaled versions such as
$ \RAWA := \CA\U-\WA\BA\U$,
where $\WA:=\MA\inv$.
We call $\MA$ the \defining{interior weight}
(or \defining{\outcrowding cap})
and we call $\WA$ the \defining{boundary weight}
(or \defining{positivity CFL cap} --- see Theorem 
\ref{positivityPreservingTheorem}).
\end{definition}
\begin{convention}[$\WA=\MA\inv$]
Throughout this paper, any symbol involving
the letter ``W'' represents the reciprocal of
the corresponding symbol obtained by replacing ``W'' with ``M''.
\end{convention}
\begin{definition}[admissible weights]
  If $\MA\ge\MAopt$, then we say that $\MA$ and $\WA$ are \defining{admissible}
  weights, because in this case positivity of $\U$ implies positivity of $\RAMA$,
  which means that enforcing positivity of $\RAMA$ respects accuracy.
\end{definition}

Formally, we have the following two theorems, which comprise
the essence of the paper.

\subsection{First fundamental theorem of outflow rate limiting}\ 

\begin{figure}[]
\fbox{\begin{minipage}{\linewidth}
\begin{center}
  \large{\textbf{I: $\h(\Um,\Up)\le\cs\Um$ for a positivity-preserving numerical flux.}}
\end{center}
\begin{center}
\includegraphics{figures/riemann.7}
\end{center}
\begin{definition}
Define the \defining{one-cell problem
with states $(Z,\Um,\Up)$, cell width $\Dx$, and flux function $f$} by 
\end{definition}
\begin{gather}
  \label{1DscalarProblem}
  \partial_t u + \partial_x f(u)=0, \qquad
 \begin{gathered}
  \u(0,x) = 
  \left\{
    \begin{array}{l l}
      Z & \hbox{ if } x<0,
   \\ \Um & \hbox{ if } 0<x<\Dx,
   \\ \Up & \hbox{ if } \Dx<x;
    \end{array}
  \right.
 \end{gathered}
\end{gather}
it is implied that $Z,\Um,\Up\ge 0$.
Define the \defining{speed cap} to be the sum of the signal speeds
entering the cell: $\cs(Z,\Um,\Up):=\sp(Z,\Um)-\sm(\Um,\Up)$.

\medskip
\begin{theorem}[outflow rate is bounded by speed cap times interior value]
\label{speedCapTheorem}
Consider the one-cell problem with states
$(Z,\Um,\Up),$ cell width $\Dx$, and positivity-preserving flux function $f$.
Define the \defining{numerical speed cap} $\cs:=\ssp(Z,\Um)-\sSm(\Um,\Up),$
where $\sSm(\Um,\Up)<\ssm(\Um,\Up)\le 0$.
Let $h_0(A,B)$ be the interface flux of the Riemann problem with
states $(A,B)$ and
let $\h(A,B)$ be a consistent numerical flux function that preserves
positivity for the Euler update
\begin{gather}
  \label{EulerUpdate}
  (\Um)\next = \Um - \tfrac{\Dt}{\Dx}\bparen{\h(\Um,\Up)-\h_0(Z,\Um)} \ge 0
\end{gather}
if $\cs\Dt\le\Dx$, i.e., if incoming signals do not cross.
Assume that $Z$, $Z\cs(Z,Z,\Um)$, and $f(Z)$ all equal or approach $0$.
Then $\h(\Um,\Up)\le\cs\Um$.
\end{theorem}
\medskip
\begin{myproof}
Choose $\Dt=\Dx/\cs$.  Then \eqref{EulerUpdate} says that
$h(\Um,\Up)\le\cs\Um+\h_0(Z,\Um)$.
But since material cannot flow out of a vacuum, $\h_0(Z,\Um)\le 0$;
to verify this rigorously, make the replacements
$\Um\to Z,$ $\Up\to\Um$, $\h\to\h_0$,
and $\cs\to\cs(Z,Z,\Um)$, and use that $\h_0(Z,Z)=f(Z)$,
which approaches $0$.
Therefore, $\h(\Um,\Up) \le \cs\Um$.
\end{myproof}

\medskip
To see how to define $\cs$ and $\h$ so that the Euler update
\eqref{EulerUpdate} maintains positivity,
for $0\le t\le\Dt$, let
$\sm_0(t)\le0\le\sp_0(t)$ be the signals emanating from $x=0$ 
and let $\sm(t)\le\Dx\le\sp(t)$
be the signals emanating from $x=\Dx$ 
(see Figure \ref{defRiemannProblem}).
For the HLL flux,
the \mention{numerical signal speeds}
$\sSm$ and $\Sp$ are required to satisfy $\sSm\le\ssm\le 0\le \ssp\le\sSp$.
Assume that $\capspd$ is an upper bound on the sum of the
signal speeds $\ssp_0$ and $-\sSm$.
Since $\Dx=\cs\Dt$,
the signals $\sp_0$ and $S^\pm(t):=\Dx+t\dot{S}^\pm$
do not cross
(i.e., $\sp_0(\Dt)\le \Sm(\Dt)\le \sm(\Dt)$).
An exact Riemann solver uses the flux of the exact solution
at the $\Dx$ interface.  The HLL solution $\uHLL(t,x)$ equals $\u(t,x)$
except in the interval
$\Signal(t):=[\Sm(t),\Sp(t)]$,
where it equals the average value
of $\u$ in $\Signal(t)$, 
$
  U^*(\Um,\Up) := \frac{f(\Um)-f(\Up)+\sSp\Up-\sSm\Um}{\sSp-\sSm}.
$
The flux needed at $x=\Dx$ to account for the amount
of material that $\uHLL$ accumulates on each side of
the interface defines the HLL numerical flux
$
  h(\Um,\Up):=\frac12\big(f(\Up)+f(\Um)+(\sSp+\sSm)\U^* -\sSp\Up-\sSm\Um\big)
$
(an instance of the formula at the bottom of Figure \ref{defRiemannProblem}).
\end{minipage}}
\caption{A signal speed cap and boundary node value bound the outgoing flux rate.
}
\label{HLLfig}
\end{figure}

The first essential theorem of outflow rate limiting is
displayed in Figure \ref{HLLfig}. It asserts that at each
boundary node the numerical outgoing flux rate $\h$ is bounded by
the product of the interior value $\Um$ at the node and a speed
cap $\cs$; it assumes that boundary node states are positive and
that $\h$ preserves positivity for a one-cell problem with
states $(0,\Um,\Up)$ and cell width $\Dx=\cs$ for any $\Dt\le 1$.
\def\Neps{\mathcal{N}_\epsilon}
\def\tNeps{\widetilde\mathcal{N}_\epsilon}
\begin{remark}
  \emph{Multiplying by cell area over cell volume, 
    $\hA(\Um,\Up)\le\csA\Um.$}
\end{remark}
\begin{remark}
\emph{Examples of such a positivity-preserving
flux function $\h$ are the numerical flux defined
by the exact Riemann solver or the Harten-Lax-van Leer
(HLL) \cite{article:HaLaLe83} or local Lax-Friedrichs (LLF)
\cite{article:Ru61} approximate Riemann solvers used with
numerical signal speeds that are truly upper bounds
for physical signal speeds. See Figure \ref{HLLfig}.
LLF is the special case of HLL where the left-going
and right-going numerical signal speeds are equal:
$\sSm+\sSp=0$.}
\end{remark}
\begin{remark}
\emph{If outflow capping is used, then it is not necessary to 
compute the speed cap $\capspd$; it is sufficient to know that
such a finite cap exists.}
\end{remark}
\begin{remark}
  \emph{The proof assumes that Riemann problems are well-defined
  for vacuum states and involve finite wave speeds. 
  Riemann problems involving vacuum states were considered
  for gas dynamics in \cite{liuSmoller80}.
  For the proof it is enough that there exist a sequence of
  strictly positive states $Z_m$ approaching vacuum such that
  the $\limsup$ of the signal speeds of the Riemann problems
  with states $(Z_m,\Um)$ is finite.
  In systems such as shallow water and gas dynamics
  for which any state can be connected to the vacuum state
  without shocks,
  one can choose states $Z_m$ so that $\ssp(Z_m,\Um)=\maxspd(\Um)$
  for all $Z_m$; 
  that is, $\ssp(0,\Um)=\maxspd(\Um)$.
  }
\end{remark}

\subsection{Second fundamental theorem of outflow rate limiting}\ 

The second essential theorem is displayed in Figure
\ref{CellOutflowRateLimiting} and states that
if the retentional $\RAWA=\CA-\WA\BA$ is positive
then an Euler step maintains positivity of the cell average
for any time step whose nondimensionalized version
$\csA\Dt$ is no greater than $\WA$, where
$\cs$ is a uniform upper bound on the speed cap
over all nodes.

\begin{figure}[]
\fbox{\begin{minipage}{\linewidth}
\begin{center}
  \large{\textbf{II: A speed cap and positive retentional guarantee a positivity-preserving time step.}}
\end{center}
  \thmskip
  \begin{theorem}
    \label{positivityPreservingTheorem}
    Suppose that:
    \begin{alignat*}{3}
      \hA&\le\csA\Um
        && \hbox{ (numerical flux is bounded via the speed cap $\csA$)}, \\
      \csA\Dt&\le\WA
        && \hbox{ (positivity CFL number is bounded by $\WA$)}, \and \\
      \WA\BA\U&\le\CA\U
        && \hbox{ (the \outcrowding is capped by $\MA=\WA\inv$).}
    \end{alignat*}
    Then the updated cell average $\CA\U-\Dt\BA\hA$ is positive.
  \end{theorem}

  \thmskip
  \begin{myproof}
    $
      \Dt\BA\hA \le \csA\Dt\BA\U \le \WA\BA\U \le \CA\U.
    $
  \end{myproof}

    \thmskip
    \emph{Remarks.} The proof says that the loss is at most the cell content.
    The retentional $\RAWA\U:=\CA\U-\WA\BA\U$
    represents the content retained if the maximum possible loss occurs
    and is positive precisely when the \outcrowding
    $\BA\U/\CA\U$ is capped by $\MA$.
    Theorem \ref{positivityPreservingTheorem} is sharp:
    a smooth example with spatially varying flux function
    satisfying $\f=\nhat\cs\u$ on $\dK$
    shows that loss can equal cell content.
    Enforcing positivity of $\RAMA$ maintains accuracy
    for scalar positivity limiting if $\MA\ge\MAopt$.
\end{minipage}}
\caption{A speed cap and boundary average cap bound the outgoing flux.}
  \label{CellOutflowRateLimiting}
\end{figure}

  The first fundamental theorem
  can be invoked to satisfy
  the first hypothesis ($\hA\le\csA\Um$) of
  the second fundamental theorem
  if at each boundary node the numerical flux function
  is defined using a 1D Riemann problem with frozen flux
  (see Figure \ref{methodOfLinesFV}).
  The proof of Theorem \ref{speedCapTheorem}
  assumes that physical solutions to
  the two auxiliary Riemann problems are well-defined and
  maintain positivity.  For simplicity, we can impose this
  assumption rather than assume that the full conservation law
  \eqref{hypsystemB} of Figure \ref{fundamentalAssumptions} maintains positivity.
  But in fact, omitting careful justification,
  it is an implication of this work that,
  with appropriate regularity and convergence assumptions, we have:

  \begin{proposition}
  The following are equivalent:
  \begin{enumerate}
   \item \label{firstCondition}
     The full system \eqref{hypsystemB} maintains positivity.
   \item \label{secondCondition}
   Any Riemann problem with positive states and frozen
     flux function maintains positivity.
   \item \label{thirdCondition}
   The one-cell Euler update \eqref{EulerUpdate} of Figure \ref{HLLfig}
      with an HLL numerical flux function maintains positivity for
      sufficiently large $\cs$.
  \end{enumerate}
  \end{proposition}
  \emph{Justification.}
  To see that \ref{firstCondition} implies \ref{secondCondition},
  assume without loss of generality that the flux function
  is frozen at $t_0=0$ and $\x_0=0$ and in the direction $\nhat.$
  The self-similar solution $u(t,x)=\uhat(x/t)$
  to the auxiliary Riemann problem with positive states $(A,B)$
  can be approximated arbitrarily well by the solution $\tu(t,\x)$
  to a multidimensional problem \eqref{conservationLaw}
  with initial data $\tu_0(\x)$ equal to $B$ if $\nhat\dotp\x>0$,
  else equal to $A$.  In particular,
  $\uhat(\xi)=\lim_{t\searrow0} \tu(t,t\xi\nhat)$.
  We omit a careful justification of this limit,
  and simply note that $f$ evaluated at $\nhat\dotp\x$
  approximates $\nhat\dotp\f$ arbitrarily well for sufficiently
  small $t$ and $\x$; here we rely on the fact that $\f$ is differentiable
  (hence continuous) and that $\tu$ by definition depends continuously on $\f$
  wherever $\tu$ is well-defined.
  Since $\tu$ satisfies positivity, so do $\uhat$ and $\u$.

  That \ref{secondCondition} implies \ref{thirdCondition}
  follows from the definition of HLL
  in terms of averagings of physical solutions
  (see Figure \ref{HLLfig}).

  That \ref{thirdCondition} implies \ref{firstCondition}
  follows from the assumption that the 
  positivity-preserving algorithm described in this work 
  converges at least in the first-order case of
  a representation space that is constant in each cell.
  We remark that condition \ref{thirdCondition}
  gives a practical means of verifying that
  a physical system maintains positivity. \myqedhere

\subsection{Direct proof for the 1D case}
See Figure \ref{1Dproof}.
\begin{figure}[]
\fbox{\begin{minipage}{\linewidth}
\begin{center}
  \large{\textbf{Comparison with Zhang and Shu (1D case)}}
\end{center}
\vskip1ex
A corollary of the theorems in Figures \ref{HLLfig} and
\ref{CellOutflowRateLimiting} is the 1D case:

\medskip
\begin{corollary}
  The Euler update \eqref{eqn:limited_update}
  maintains positivity of the cell average if the retentional
  $\RA\U^n_i:=\Uavg^n_i-\WA\BA\U^n_i$ is positive in each cell
  and if the time step satisfies $\frac{2\Dt\cs}{\Dx}\le\WA$,
  where $\cs/2$ is an upper bound on signal speeds and
  $\BA\U_i=\frac{\Ucm+\Ucp}2$ is the boundary average;
  we assume that the frozen flux function $\hp$ used at
  any interface is positivity-preserving,
  in the sense that if it is used at all interfaces then it maintains
  positivity of cell averages for data constant in each cell if
  $\cs\frac{\Dt}{\Dx}\le1$
  (i.e., $\Dt$ is short enough that signals from cell interfaces
  cannot cross in Godunov's method).
\end{corollary}

\medskip
To facilitate comparison with
Theorem 2.1 of \cite{article:ZhangShu10rectmesh},
we offer a direct proof.

\medskip
\begin{myproof}
Substituting $\Uavg^n_i= \RA\U_i + \WA\BA\U_i$
(a boundary-weighted quadrature rule),
\begin{align*}
   \Uavg\next_i &= \Uavg^n_i - \tfrac{\Dt}{\Dx}\bparen{\hp(\Ump,\Upp)-\hm(\Umm,\Ump)}
     \\
     &\ge \RA\U_i + \frac{\WA}{2}\bparen{
         \begin{aligned}
           &\Ucp - \tfrac{2\Dt}{\WA\Dx}\bparen{\hp(\Ucp,\Upp)-\hp(\Zcp,\Ucp)}
           \\ +                 
           &\Ucm - \tfrac{2\Dt}{\WA\Dx}\bparen{\hm(\Ucm,\Zcm)-\hm(\Umm,\Ucm)}
         \end{aligned}
        },
\end{align*}
which is a positive combination of positive quantities if
the Euler steps in brackets satisfy 
$\frac{2\cs\Dt}{\WA\Dx}\le 1$;
here we take $\Zcp$ and $\Zcm$ to be the vacuum state
and we use that material cannot flow out of a vacuum:
$\hp(\Zcp,\Ucp)\le 0$ and
 $\hm(\Ucm,\Zcm)\ge 0$.
Zhang and Shu instead assume a spatially invariant flux
function ($\hp=\hm$) and choose $\Zcp=\Upm$ and $\Zcm=\Ump$.
They write the retentional as a weighted sum over interior positivity points
(see Theorem \ref{optimalBoundaryWeightedQuadratureExists})
and choose $\WA=\WAopt$ (see Theorem \ref{weightsFor1Dintervals}).
\end{myproof}

\end{minipage}}
\caption{Direct proof for the 1D case.}
\label{1Dproof}
\end{figure}

\subsection{Affine-invariant definitions of wave speed and retentional}
\label{secAffineInvariants}


We have used bars in Sections
\ref{sec:dg_schemes}--\ref{sec:framework} to indicate the
scale-invariant formulation. Note that all definitions and
statements remain valid if bars are dropped. Alternatively,
we can redefine barred quantities using affine-invariant
definitions.

Heretofore we have made no distinction between physical and
canonical mesh cell coordinates. In general, each physical mesh
cell is the image of a \mention{canonical mesh cell} under a
coordinate map. In the case of isoparametric mesh cells, the
coordinate map is a diffeomorphism, and it is important to work
in canonical coordinates.  If the coordinate map is affine,
however, then we can avoid making a distinction between physical
and canonical coordinates if we adopt an \mention{affine-invariant}
formulation.

In canonical coordinates, the mesh cell is almost always a cube
or simplex and the representation space almost always consists
of polynomials.  This remains true for isoparametric mesh cells.
With an affine-invariant formulation, our results are generally
applicable, independent of whether the canonical simplex is
defined to be regular or the corner of a box.

In Section \ref{CalculationOfWeights}, we tabulate values or estimates
of $\MAopt$ for \emph{regular} canonical polytopes.
Thus, when using our tabulated weights to apply the outflow
positivity limiting framework, one should take wave
speeds to be in terms of their values in the coordinates of a regular
canonical mesh cell. A natural way to do this is to use affine
invariants.

An affine-invariant formulation of the Euler step \eqref{explicitEuler}
that agrees with the scale-invariant formulation \eqref{explicitEulerA}
in the case of a regular canonical mesh cell is
$ 
  \CA\U\next = \CA\U - \Dt\BAF\hAF,
$ 
where \emph{$\BAF$ 
is an affine-invariant version of the boundary average
that computes the arithmetic average 
over all faces 
of the average on each face};
here
$ \hAF_\e:=\frac{|\faces|\dA_\e}{\vol}\h_\e $
is an affine-invariant version of the numerical flux,
where $\dA_\e$ is the area of the face of boundary node $\e$
and $|\faces|$ is the number of faces of the mesh cell.
Similarly, since wave speeds scale like fluxes,
if $\cs_\e$ is a wave speed at boundary node $\e$ in
the canonical coordinates of mesh cell $\K$, then
the quantity 
$ 
  \csAF_\e:=\frac{|\faces|\dA_\e}{\vol}\cs_\e
$ 
is the affine-invariant version
that agrees with the scale-invariant version
$\csA_\e=\frac{\area}{\vol}\cs_\e$ in the case of a canonical regular polytope.
In the case of a regular canonical polytope
the affine-invariant quantities
$\hAF_\e$, $\csAF_\e$, and $\BAF$
agree with their scale invariant equivalents
$\hA_\e$, $\csA_\e$, and $\BA$,
and 
the theory of Sections \ref{sec:framework}--\ref{sec:accuracy}
goes through.

If one is computing with a non-regular canonical polytope $\tK$
(e.g.\ the standard orthogonal simplex defined as the corner of
a box) then an affine-invariant definition should be used in the
definition of the retentional:
\begin{equation}
  \label{affineInvariantDefs}
 \begin{aligned}
  \RWA(\U) &:= \CA\U - \WA\BAF\U.
 \end{aligned}
\end{equation}
If outflow capping is used, then computing wave speeds is not
necessary, and use of the boundary face average $\BAF$ in 
the definition \eqref{affineInvariantDefs} of the retentional
is the only modification needed to implement
affine-invariant positivity-limiting.

\section{Systems versus scalar case}
 \label{sec:systems}

In the systems case of Figure \ref{fundamentalAssumptions},
the set of positive states $\calP$ is a convex set
$\calP\subset\realsN$.
A corollary of Theorem 3.4 in \cite{rudinFA73} is that
any open or closed convex set is an intersection of half-planes.
But any half-plane is the set on which some affine functional $\Alpha$ is positive.
Therefore, we assume that there exists a set $\dP$ of affine functionals
such that a state $\vu\in\realsN$ is positive if for all $\Alpha\in\dP$
$\Alpha(\vu)\ge 0$.  We call $\Alpha$ a \defining{state positivity functional}.
Any such affine functional $\Alpha$ decomposes uniquely as
$\Alpha=:s+\Lambda$, where $s=\Alpha(0)$ is a scalar shift
and $\Lambda:=\Alpha - s$ is its \defining{linear part}.
Applying $\Alpha$ reduces the systems case to the scalar case:
\begin{align*}
  \partial_t \Alpha(\vu) + \Div(\Lambda\vbf)=0.
\end{align*}
After applying $\Alpha$ to states and $\Lambda$ to fluxes,
all statements and reasoning of Theorems
\ref{speedCapTheorem}
and \ref{positivityPreservingTheorem}
remain valid.
For example, the inequality $\h\le\cs\Um$ becomes
the inequality $(\Lambda\vh)\le\cs\Alpha(\vUm).$

Let $\R:\Uspace\to\reals$ be a retentional.
Applied component-wise, the retentional
defines a state-valued linear map
$\tR:\Uspace\to\realsN$, which we identify with $\R$.
Applied point-wise, a state positivity functional
$\Alpha=s+\Lambda$
defines a map $\tAlpha=s+\tLambda:\Uspace^N\to\Uspace$
that we identify with $\Alpha$.
Observe that
$\R\tLambda=\Lambda\tR$.
It is desirable that
the retentional commute with all state positivity functionals:
$\R\tAlpha=\Alpha\tR$;
then enforcing positivity of $\tR$ is the same as enforcing
$\R\tAlpha\ge0$ for all $\Alpha\in\dP$.
This will hold if $\R s=s,$
which holds if $s$ is always $0$
(that is, if $\calP$ is a convex \emph{cone})
or if the retentional is rescaled so that
$\R 1= 1$.  Therefore, for any retentional $\R$
we define its \defining{unital retentional}
to be $\Rhat:=\frac{\R}{\R 1}$.
Then $\Rhat 1=1$.  
To be concrete:
$\Rhat^M=\frac{\M\C-\B}{\M-1}$;
we remark that $\Rhat$ is an affine combination of state values
and therefore is invariant under translation of state space.

For typical systems such as shallow water and gas dynamics,
the set of positive states is a convex cone (i.e.\ invariant
under rescaling), so $s=0$ and $\Alpha=\Lambda$ and there is no need
to rescale the retentional to its unital version.
We remark that by adding the trivial equation
$\partial_t\uextra=0$ to the system, where $\uextra$
is a scalar taken to be $1$ for physical solutions,
defining the set of positive states to be
$\setof{(r\u,r):\u\in\calP,r\ge0}$,
and extending each positivity functional $\Alpha=s+\Lambda$
to the linear functional $\tLambda:=(\u,\uextra)\mapsto s\uextra + \Lambda\u$,
the set of positive states can be assumed without loss of
generality to be a convex cone.

While we can assume the scalar case without loss of generality
when it comes to the \emph{positivity-preserving} results of
Section \ref{sec:framework},
the \emph{accuracy} result of Theorem \ref{OptimalWeightPreservesAccuracy}
does not go through e.g.\ to the gas-dynamics case
(see Section \ref{AccuracyOfPositivityLimitingInTheSystemsCase}).


\begin{figure}[]
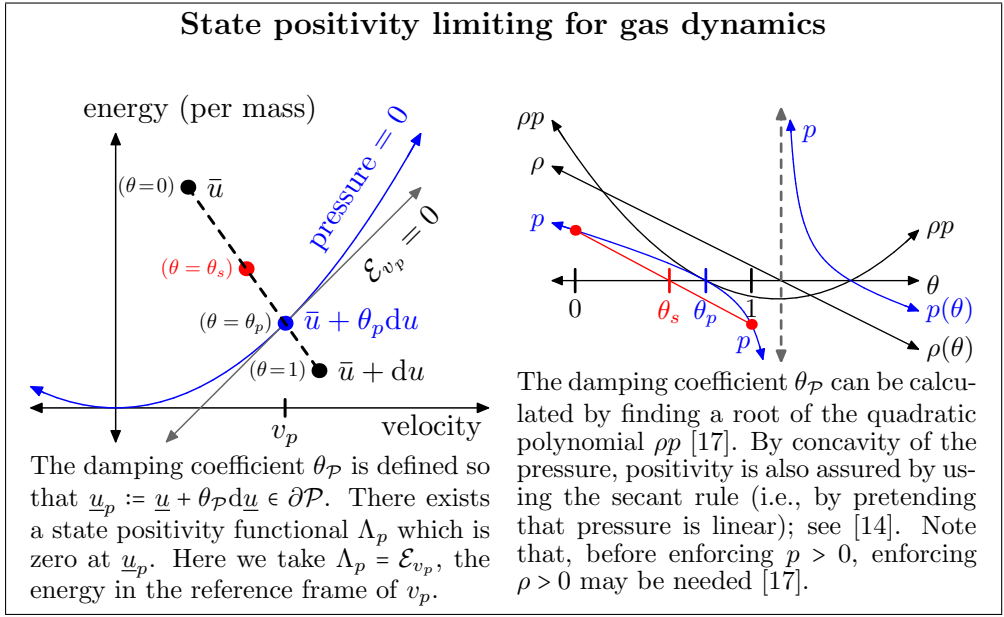

\fbox{\begin{minipage}{\linewidth}
\begin{center}
  \large{\textbf{State positivity limiting for gas dynamics}}
\end{center}
\begin{center}
\begin{gather*}
 \begin{matrix}
   \begin{matrix}
    \begin{matrix}
      \mbox{\includegraphics[width=.47\linewidth]{figures/limiting.24}}
    \end{matrix}
    \\
    \begin{matrix}
      \mbox{\parbox{.47\textwidth}{
       The damping coefficient $\thetaP$ is defined
       so that $\vu_p:=\vu+\thetaP\vdu\in\partial\calP$.
       There exists a state positivity functional $\Lambda_p$
       which is zero at $\vu_p$. Here we take
       $\Lambda_p=\energy_{v_p}$, the energy in the reference
       frame of $v_p$.
       }}
    \end{matrix}
   \end{matrix}
\quad
   \begin{matrix}
    \begin{matrix}
      \mbox{\includegraphics[width=.47\linewidth]{figures/limiting.11}}
    \end{matrix}
    \\
    \begin{matrix}
      \mbox{\parbox{.47\textwidth}{
       The damping coefficient $\thetaP$ can be calculated by finding
       a root of the quadratic polynomial $\rho p$
       \cite{article:ZhangShu10rectmesh}.
       By concavity of the pressure, positivity is also assured by
       using the secant rule (i.e., by pretending that pressure is linear);
       see \cite{WangZhangShuNing12}.
       Note that, before enforcing $p>0$,
       enforcing $\rho>0$ may be needed
       \cite{article:ZhangShu10rectmesh}.
       }}
    \end{matrix}
   \end{matrix}
 \end{matrix}
\end{gather*}
\end{center}
\end{minipage}}
 \caption{Depiction of positivity limiting of a state value
   for Euler gas dynamics.  Positivity limiting requires determining
   where a linear path in state space intersects the boundary of positive states.}
 \label{positivityLimitingForSystems}
\end{figure}

\subsection{State positivity indicators}\ 
\label{StatePositivityIndicators}

Outflow positivity requires the solution to the following problems:
\begin{enumerate}
\item To cap outflow,
  determine the largest value of $\theta\in[0,1]$
  for which $\u(\theta):=\ubar+\theta\du$ satisfies positivity,
  where $\ubar=\C\U$ and $\du=\Dtmax\B\h$.
\item To cap \outcrowding at $\MA$,
  determine the largest value of $\theta\in[0,1]$
  for which $\u(\theta):=\ubar+\theta\du$ satisfies positivity,
  where $\ubar=\Rhat\Uavg$, $\du=\Rhat\dU$, 
  $\Rhat:=(\RAMA 1)\inv\RAMA$, and $\RAMA:=\MA\CA-\BA$.
\item To enforce positivity at $\x_0$,
  determine the largest value of $\theta\in[0,1]$
  for which $\u(\theta):=\ubar+\theta\du$ satisfies positivity,
  where $\ubar=\Uavg$ and $\du=\dU(\x_0)$.
\end{enumerate}
All these problems seek the intersection of a line
segment with the boundary of a convex set
of states designated positive and
require computing a damping coefficient:
\begin{definition}
  Let $\CP$ be a closed convex set of states deemed positive,
  and let $\dP$ be a collection of affine functionals such that
  $\CP=\setof{\u\in\realsN: (\forall\Alpha\in\dP)\ \Alpha\u\ge 0}$.
  Let $\Alpha\in\dP$.
  Let $\ubar,\u\in\realsN$, where $\ubar\in\CP$.
  Define $\du:=\u-\ubar$.
  Define $\thetaA(\vu,\vdu)$ to be the largest $\theta\in[0,1]$
  such that $\Alpha(\thetabar\ubar+\theta\u)\ge 0$,
  where we define $\thetabar:=1-\theta$.
  Define $\thetaP(\ubar,\u)$ to be the largest $\theta\in[0,1]$
  such that $\thetabar\ubar+\theta\u\in\CP$.
  Observe that $\thetaA(\ubar,\u)$ equals $1$ if $\Alpha(\u)\ge0$
  and else equals 
  $\frac{\Alpha\ubar}{\Alpha\ubar-\Alpha\u}=\frac{s+\Lambda\ubar}{-\Lambda\du}$.
  Observe that $\thetaP(\vu,\vdu)=\min_{\Alpha\in\dP}\thetaA(\vu,\vdu)$.
\end{definition}

If $\dP$ is a finite set, then it can be used directly to calculate $\thetaP$,
as illustrated in Figure \ref{outflowCappingInTheScalarCase}
for the case of scalar outflow capping.
%
In the example of Euler gas dynamics, however, the state is positive
if the energy in an arbitrary reference frame (a linear
functional of conserved variables) is positive. But such energies
comprise an infinite collection $\dP$ of linear functionals. Instead,
one can use the density $\rho$ (a linear functional), and
pressure $p$ (a concave nonlinear functional, see section 3.1 of
\cite{WangZhangShuNing12}) or $\rho p$ (a quadratic functional
of conserved variables,
see \cite{article:ZhangShu10rectmesh}) as a finite set of
positivity indicators.
See Figure \ref{positivityLimitingForSystems}.

\subsection{Accuracy of positivity limiting in the systems case}
\label{AccuracyOfPositivityLimitingInTheSystemsCase}


For the case of scalar conservation laws,
accuracy of pointwise positivity limiting
via linear damping is established in Section
\ref{sec:accuracy} based on Theorem \ref{OptimalWeightPreservesAccuracy}.

In the case of shallow water, positive states are characterized
by positivity of a single linear functional (the depth).
Therefore, the proof of accuracy for the scalar case can be
invoked to conclude that positivity limiting does not compromise
accuracy of the water depth. Assuming that fluid velocities are
bounded, this in turn implies that positivity limiting does not
compromise accuracy of the solution as a whole.

But in the case of gas dynamics, one must also enforce positivity
of the pressure.  There exist physical solutions with arbitrarily
large and rapid variation in density for which the
pressure remains arbitrarily close to zero.  Therefore,
if positivity is enforced by damping the deviation of the cell average by a
scalar value of $\theta\in[0,1]$, then 
a given amount of damping needed to enforce positivity of the pressure
can entail damping of the density variation that is arbitrarily large in
magnitude.

On the other hand, if the pressure of the exact solution
is strictly bounded away from zero, then positivity limiting
respects accuracy for the simple reason that positivity
limiters are not triggered if the mesh is sufficiently fine,
if the exact solution is smooth, and if $\MA>1$.
In this case, positivity limiting
is really about practical robustness.
%
The study of accuracy and stability for gas-dynamic problems
(or sequences of problems) for which the pressure is not
bounded away from zero is delicate and is left to 
future work.

%

\def\tU{\widetilde{\U}}
\def\uT{\u_T}
\def\uR{\u_R}
\section{Admissible weights and cell positivity functionals}
\label{sec:accuracy}


In this section we show
that positivity of the retentional $\RM(\U)$ can be
enforced without loss of accuracy as long as the
interior weight $\M$ 
is greater than or equal to an optimal value $\Mopt[\K]$
and that this holds precisely when positivity of $\RM(\U)$
can be enforced by enforcing positivity at appropriately chosen
\mention{interior points} in the mesh cell $\K$
(see Figure \ref{pointwiseVersusRetentionalPositivityLimiting}).

In Section \ref{CalculationOfWeights} we use this theory
(and interior points in particular)
to compute (or isolate) the optimal value $\MAopt$ for
important mesh cell geometries and representation spaces.
As elsewhere,
in this section the following context is assumed:
\begin{context}
  \label{accuracyContext}
  $\K$ is a compact mesh cell (in canonical coordinates) and
  $\Uspace$ is a finite-dimensional polynomial space which
  contains the set $\PkD$ of all polynomials in $D$ variables
  of degree at most $k\ge0$.
  We always assume that $\U\in\Uspace$.

  Recall from equation \eqref{MAoptDef} that the optimal interior weight
  $\Mopt=\Mopt[\K](\Uspace)$ is defined as the supremum of the
  \outcrowding $\Bhat(\U):=\frac{\B(\U)}{\C(\U)}$ over
  the set $\Uplus[\K]$ of all nonzero solutions $\U\in\Uspace$
  positive in the mesh cell $\K$, where $\C(\U):=\intK\U$ and
  $\B(\U):=\intdK\U$. The retentional $\RW(\U):=\C(\U)-\W\B(\U)$,
  or $\RM(\U):=\M\C(\U)-\B(\U)$, is positive for all such positive
  $\U$ precisely when $\M:=\W\inv$ is \mention{admissible}, i.e.,
  when $\M\ge\Mopt$.
\end{context}

The following results will be proved in detail in Part II
\cite{article:Jo12}.
We summarize these results in Figure \ref{mutualConfirmationFig}
and illustrate them in Figure
\ref{pointwiseVersusRetentionalPositivityLimiting}.

\begin{theorem}[{An admissible weight preserves accuracy}]
\label{OptimalWeightPreservesAccuracy}
  Let $\M\ge\Mopt$.  Then linearly damping the deviation
  from the cell average
  just enough to enforce positivity of the retentional
  $\RM(\U)$ retains order-$k$ local accuracy.
\end{theorem}
\begin{proofsketch}
  Since physical solutions are positive, damping the variation
  from the cell average just enough to enforce positivity
  everywhere in the mesh cell preserves accuracy; a formal
  proof uses that the maximum $\nmax:=\dU\mapsto\max(\dU(\K))$ and
  the magnitude of the minimum $\nmin:=\dU\mapsto\max(-\dU(\K))$
  are (equivalent) asymmetric norms on the finite-dimensional linear space of
  polynomials $\dU\in\Uspace$ whose cell average is zero
  (see \cite{article:ZhangShu11}).
  If $\M$ is an admissible weight, then
  the retentional $\RM(\U)$ is positive if $\U$ is positive
  everywhere in the mesh cell,
  so enforcing positivity of $\RM(\U)$ also retains accuracy.
\end{proofsketch}

\begin{theorem}[{Enforcing positivity of the
  points of a nodal scheme is sufficient to enforce
  a positive retentional for some $\M<\infty$}]
\label{existenceOfFiniteM}
  Let $\X$ be a set of points capable of representing 
  the solution and for which the cell average can be represented
  as a strictly positive combination of point values.
  Then $\MX<\infty$, where
  $\MX:=\sup\Bhat(\Uplushat[\X])$ and
  $\Uplushat[\X]:=\setof{\U\in\Uspace: \C\U=1 \and (\forall \x\in\X)\ \U(\x)\ge0}$.
\end{theorem}
\begin{proofsketch}
  If $\U$ is positive on $\X$
  then $\C\U>0$, so $\Bhat(\U)$ is continuous on $\Uplushat[\X]$.
  Using that $\nmax$ and $\nmin$ are equivalent asymmetric norms,
  one can show that $\Uplushat[\K]:=\setof{\U\in\Uplus[\K]:\C\U=1}$
  is bounded and hence compact.
  Therefore, $\Bhat(\Uplushat[\X])$ has a finite maximum.
\end{proofsketch}

\begin{theorem}[{A boundary-weighted quadrature rule gives an
  upper bound for the optimal weight}]
  Suppose that the cell content
  is given by a boundary-weighted quadrature rule,
  $\C(\U) = \R(\U) + \W\B(\U)$, where the boundary weight $\W$
  is strictly positive and $\R(\U)$ represents a quadrature rule
  $\sum_{\x\in\X} \wx \U(\x)$ for the retentional $\RW(\U)$ with positive weights $\wx$
  that is exact for $\U\in\Uspace$.
  Then $\Mopt\le\M:=\W\inv$.
\label{upperBoundViaQuadratureRule}
\end{theorem}
\begin{remark}
  \label{retentionalPositivityPoints}
\emph{Enforcing positivity at the points $\X$ enforces positivity of
  $\RW$; we call $\X$ a set of
  \defining{retentional} (or \defining{interior}) \defining{points}
  for the weight $\M$.}
\end{remark}
\begin{myproof}
  Indeed, $\R=\C-\W\B$ is positive for all $\U$ positive in $\K$,
  so for any positive $\U$,
  $\W\inv\ge\frac{\B(\U)}{\C(\U)}=\Bhat(\U)$,
  so $\W\inv$ is an upper bound on $\Mopt$.
\end{myproof}

\begin{theorem}[{An optimal boundary-weighted quadrature rule exists}]
  \label{optimalBoundaryWeightedQuadratureExists}
  An optimal quadrature rule
  $\C(\U)=\sum_x \wx\U(\x) + \Wopt\B(\U)$ exists
  such that $\Mopt=\Wopt\inv$.
\end{theorem}
\begin{proofsketch}
  If $\X$ is a subset of $\K$, then
  $\RM$ is positive on $\Uplushat[\X]$
  if $\M\ge\Mopt[\X]:=\sup\Bhat(\Uplushat[\X])$ and
  thus (if $\X$ is a finite set) is representable as a positive
  combination of values at a subset $\tX\subset\X$ of points
  which comprise linearly independent point evaluation functionals;
  a compactness argument extends this statement from finite $\X$
  to $\X=\K$.
\end{proofsketch}

\begin{theorem}
  \label{optimizerExists}
  An optimizer $\Uopt\in\Uplus[\K]$ exists such that
  $\Mopt=\Bhat(\Uopt)$.
\end{theorem}
\begin{proofsketch}
  By definition, $\Mopt$ is defined to be the supremum of the
  \outcrowding
  $\Bhat(\U):=\frac{\B(\U)}{\C(\U)}$
  over all nonzero $\U$ positive in the mesh cell $\K$.
  Thus, any $\U$ positive on $\K$ gives a lower bound $\Bhat(\U)$ on $\Mopt$.
  The question is whether an \argmax of $\Bhat$ exists.
  %
  Since $\Bhat$ is continuous $\Uplushat[\K]$
  and $\Uplushat[\K]$ is compact
  (see the proof of Theorem \ref{OptimalWeightPreservesAccuracy}),
  an \argmax $\Uopt$ of $\Bhat$ exists.
\end{proofsketch}


Thus, it is possible to determine the optimal weight $\Mopt$
and a set of optimal retentional points
simply by guessing an optimizer and an optimal boundary-weighted
quadrature rule and using them to confirm one another.
We now identify properties that reduce the set of guesses that
we have to consider.

\begin{theorem}[{Invariance under isometries can be required
of retentional points and optimizer candidates}]
\label{invarianceUnderIsometryCanBeRequired}
\end{theorem}
\begin{proofsketch}
A typical canonical mesh cell has a group of isometries;
by averaging over an orbit, for any solution $\U$ there exists
a $\tU$ invariant under the action of the group of isometries of
the mesh cell with the same values of $\C$, $\B$, and $\Bhat$.
Similar averaging shows that every boundary-weighted quadrature rule
  $\C(\U)=\sum_{\x\in\X} \wx\U(\x) + \W\B(\U)$
has an invariant version
  $\C(\U)=\sum_{\x\in\tX} \twx\U(\x) + \W\B(\U)$.
\end{proofsketch}

\def\ta{s}
\def\tb{t}
\def\Ua{U}
\def\Ub{V}
\begin{lemma}[The \outcrowding of a positive combination
 of two functions with positive cell average is a convex
 combination of their \outcrowding values]
  \label{BhatOfSumIsConvexCombination}
  Let $\Ua,\Ub\in\Uplus[\K]$.
  Let $\ta>0$ and $\tb>0$.
  Then $\Bhat(\ta\Ua+\tb\Ub) = a\Bhat(\Ua)+(1-a)\Bhat(\Ub)$
  for some $a\in(0,1)$.
\end{lemma}
\begin{myproof}
  $\Bhat(\ta\Ua+\tb\Ub)
    =\frac{\ta\B\Ua+\tb\B\Ub}
          {\ta\C\Ua+\tb\C\Ub}
    = a\Bhat(\Ua)+(1-a)\Bhat(\Ub),
  $
where $a = \frac{\ta\C\Ua}
                {\ta\C\Ua+\tb\C\Ub}. 
      $
\end{myproof}

\begin{theorem}[A non-constant optimizer has a zero]
  \label{nonEmptyZeroSet}
  An optimizer $\Uopt$ has a zero in $\K$ unless
  $\U=1$ is an optimizer (in which case $\MAopt=1$).
\end{theorem}
\begin{myproof}
  We have that $\MAopt = \BAhat\Uopt.$
  Let $\Umin:=\min\Uopt(\K)$.
  Assume that $\Umin>0$.
  Let $\tU:= \Uopt-\Umin$.
  Since $\BAhat(\Umin)=1$,
  by Lemma \ref{BhatOfSumIsConvexCombination},
  $\BAhat(\Uopt) = \BAhat(\tU+\Umin)
    = a\BAhat(\tU) + (1-a)$,
  for some $0<a<1$.
  If $\BAhat(\tU)<1$
  then $\BAhat(\Uopt)<1=\BAhat(\Umin)$, contradicting
  the maximality of $\BAhat(\Uopt)$.
  If $\BAhat(\tU)>1$
  then $\BAhat(\Uopt)<\BAhat(\tU)$, again contradicting
  the maximality of $\BAhat(\Uopt)$.
  Therefore, $\BAhat(\Uopt)=1$.
  But $\BAhat(1)=1$, so $1$ is an optimizer.

\begin{figure}[]
\fbox{\begin{minipage}{\linewidth}
\begin{center}
  \large{\textbf{Optimal accuracy-respecting cell positivity functionals}}
\end{center}
\vskip1ex

Assume Context \ref{accuracyContext}:
$\K$ is a compact mesh cell,
$\Uspace$ is a finite-dimensional polynomial space which
contains the set $\PkD$ of all polynomials in $D$ variables
of degree at most $k\ge0$,
$\CA\U:=\avgK\U,\,$
$\BA\U:=\avgdK\U,\,$
$\BAhat\U:=\CA\U/\BA\U,\,$
$\Uplus[\K]:=\setof{\U\in\Uspace: \inf\U(\K)\ge0,\,\CA\U>0},\,$
$\MAopt:=\sup\BAhat(\Uplus[\K]),$ and
$\RMA:=\MA\CA-\BA$.  Then the following results hold:

\thmskip
\begin{corollary}[\textbf{existence of an invariant
  optimizer and of invariant retentional points
  which isolate the optimal weight}]
  Linearly damping the deviation from the cell average
  just enough to enforce positivity of the retentional
  $\RMA$ retains order-$k$ local accuracy as long as
  $\MA$ is admissible, i.e., $\MA\ge\MAopt$
  (Theorem \ref{OptimalWeightPreservesAccuracy}).
  Any solution $\U\in\Uspace$ positive in $\K$ and
  any correct boundary-weighted quadrature rule
  $\CA(\U)=\sum_{\x\in\X} \wAx\U(\x) + \WA\BA(\U)$ give a bracket
  $\BAhat(\U)\le\MAopt\le\WA\inv$ for the optimal weight
  (Theorem \ref{upperBoundViaQuadratureRule}).
  Furthermore, an optimizer $\Uopt$ exists
  (Theorem \ref{optimizerExists})
  and an optimal boundary-weighted
  quadrature rule exists
  (Theorem \ref{optimalBoundaryWeightedQuadratureExists})
  such that the bracket is an equality:
  $\BAhat(\Uopt)=\MAopt=\WA\inv$.
  Enforcing positivity at the set $\X$ of
  \defining{retentional points}
  enforces positivity of
  the retentional $\RMA$
  (see Remark \ref{retentionalPositivityPoints}).
  The zero-set of $\Uopt$ is nonempty if $\MAopt\ne1$
  (Theorem \ref{nonEmptyZeroSet})
  and must contain $\X$ if $\X$ is optimal
  (by Theorem \ref{optimalBoundaryWeightedQuadratureExists}).
  All these statements continue to hold if
  the optimizer $\Uopt$ or the quadrature rule
  (i.e.\ the retentional points $\X$
  and the quadrature weights $\wAx$)
  is required to be invariant under the action of the
  isometries of the mesh cell (Theorem
  \ref{invarianceUnderIsometryCanBeRequired}).
  \label{mutualConfirmation}
\end{corollary}
\end{minipage}}
 \caption{Summary of the results of Section \ref{sec:accuracy}.}
 \label{mutualConfirmationFig}
\end{figure}

\end{myproof}

%


\section{Calculated weights (functional analysis calculations)}
\label{CalculationOfWeights}

The previous section developed a general theory of cell
positivity functionals.
%
In this section, we use the results of Section \ref{sec:accuracy}
summarized in Figure \ref{mutualConfirmationFig}
to determine retentional points and
lower and upper bounds for the optimal weight
$\MAopt[\K]$ for the two standard mesh cell geometries
(box and simplex) and for polynomial representation
spaces of varying order.
The most important results of this section are summarized
in Figure \ref{tabulatedWeights}.

\begin{remark}[retaining bars when computing weights]
\emph{In accordance with Remark \ref{droppingBars},
we dropped bars in the previous section without
affecting the validity of the general theorems.
In this section, however, we obtain concrete numerical
values for the interior weight $\MA$; since we want these results
to be independent of the scale of the canonical mesh cell, we
retain bars in this section.}
\end{remark}

\def\WAsphDk{\WAsph{D}{k}}
\def\WAinterval{\WA_{\unitinterval}}
\def\WAintervalk{\WA^k_{\unitinterval}}
\newcommand\WAbox[2]{\ensuremath{\WA^{#2}_{\unitinterval^{#1}}}}
\newcommand\WAboxD[1]{\WAbox{#1}}
\def\WAboxDk{\WAbox{D}{k}}
\newcommand\WAtri[2]{\ensuremath{{\WA^{#2}_{\medtriangleup^{#1}}}}}
\def\WAtriDk{\WAtri{D}{k}}
\newcommand\WAsph[2]{\ensuremath{{\WA^{#2}_{\bigcirc^{#1}}}}}
\def\mysph#1{\bigcirc^{#1}}
\def\mybox#1{\unitinterval^{#1}}
\def\mytri#1{\bigtriangleup^{#1}}
\def\ostar#1{\largestar^{#1}} 
\newcommand\Usph[2]{\ensuremath{{\U_{\mysph{#1}}^{#2}}}}
\newcommand\Ubox[2]{\ensuremath{{\U_{\mybox{#1}}^{#2}}}}
\newcommand\Utri[2]{\ensuremath{{\U_{\mytri{#1}}^{#2}}}}

\newcommand\Uoptsph[2]{\ensuremath{\Usph{#1}{#2,\star}}}
\newcommand\Uoptbox[2]{\ensuremath{\Ubox{#1}{#2,\star}}}
\newcommand\Uopttri[2]{\ensuremath{\Utri{#1}{#2,\star}}}
\newcommand\UoptB[2]{\ensuremath{\U_{#1}^{#2,\star}}}

\def\MAsphDk{\MAsph{D}{k}}
\def\MAboxDk{\MAbox{D}{k}}
\def\MAtriDk{\MAtri{D}{k}}
\subsection{Definitions and summary of results}

Recall that $\MAopt[\K](\Uspace)$ denotes the optimal interior
weight for mesh cell $\K$ with representation space $\Uspace$.
We tabulate values or upper and lower bounds
on the optimal interior weights
for standard canonical mesh cells:
the unit interval $[0,1]$,
the unit square $[0,1]^2$,
the unit cube $[0,1]^3$,
an equilateral triangle $\bigtriangleup^2$,
and a regular tetrahedron $\bigtriangleup^3$.
Each upper bound $\MA$ constitutes an admissible weight
and is calculating using a quadrature rule for the retentional
$\RMA$ whose quadrature points can be used as positivity
points.

\begin{definition}[canonical mesh cells]
\emph{In $\reals^D$, define canonical mesh cells:
   \begin{gather*}
     \begin{aligned}
       \mybox{D}&= &&\textrm{regular box}, &\qquad
       \mysph{D}&= &&\textrm{sphere}, \\
       \mytri{D}&= &&\textrm{regular simplex}, &\qquad
       \ostar{D}&: &&\textrm{``star-regular'' polytope}.
     \end{aligned}
   \end{gather*}
  We define a \defining{star-regular} polytope
  to be any star-convex mesh cell that 
  can be centered on the origin to have a constant value of
  $\nhat\dotp\x$, where $\nhat$ is the outward
  unit normal and $\x$ is position on the cell boundary.
  We use $\ostar{D}$ as a generic designation for
  a star-regular polytope.  Regular polytopes are star-regular.}
\end{definition}
Since most often the representation space is $\PkD$, we define 
\begin{gather*}
\MA{}_{\K}^k := \MAopt[\K](\PkD).
\end{gather*}

\begin{figure}[]
\fbox{\begin{minipage}{\linewidth}
\begin{center}
  \large{\textbf{Tabulated results of Section \ref{CalculationOfWeights}}}
\end{center}

 \subsubsection{Boxes and even star-regular mesh cells}

 If the mesh cell $K$ is even (i.e., there exists $\Kbar\in\K$ such that
 $\K\subset\Kbar-\K$), then for a regular
 polygon the optimal interior
 weight for polynomials of degree at most $k$
 is bounded above by the optimal interior weight
 $\MAsph{D}{k}$ for a $D$-dimensional sphere.
 Based on equations \eqref{WAsph} and \eqref{WAboxplusRecurrence},
 we tabulate exact values or bounding intervals for
 the optimal interior weight for an interval ($\MAintervalk$),
 a square ($\MAbox{2}{k}$), and a cube ($\MAbox{3}{k}$).
 Note that for a tensor product polynomial space
 the optimal weight is $\MAintervalk$ independent of the
 dimension $D$ of the box (see Remark \ref{tensorProductWeights}).
 Note also that for $k\le 3$, \emph{enforcing positivity at the cell
 center enforces positivity of the retentional for the optimal weight}
 (see Theorems \ref{quadraticRepresentationSpace}
 and \ref{oddDegreeForBoxes}).

   \begin{equation*}
   \begin{array}{l||c|c|c|c|c|c|}
          k\ (\Uspace=\PkD):
          & 0,1
          & 2,3
          & 4,5
          & 6,7
          & 8,9
          & 10,11
          \rule{.0pt}{2.5ex}
     \\ \hline
           n=\lfloor k/2\rfloor:
           & 0 & 1 & 2 & 3
           & 4 & 5
     \\ \hline
           m=\lfloor n/2\rfloor:
           & 0 & 0 & 1 & 1
           & 2 & 2
     \\ \hline
        \hline
        \MAintervalk = \frac{(n+1)(n+2)}{2}:
          & 1
          & 3 
          & 6 
          & 10
          & 15
          & 21
          \rule{.0pt}{2.5ex}
     \\ \hline
        \MAsph{2}{k} = (m+1)\cdot\lfloor\frac{n+3}{2}\rfloor
          & 1
          & 2 
          & 4 
          & 6
          & 9
          & 12
          \rule{.0pt}{2.5ex}
     \\ \hline
        \MAsph{3}{k} = \frac{m+1}{3}\cdot\left(3+2\lfloor\frac{n+1}{2}\rfloor\right)
          & 1
          & 1.\overbar{6} 
          & 3.\overbar{3} 
          & 4.\overbar{6} 
          & 7
          & 9
          \rule{.0pt}{2.5ex}
     \\ \hline
        \hline
        \MAbox{2}{k}
          & 1
          & 2 
          & [3.5,4]
          & [5.5,6]
          & [8,9]
          & [11,12]
          \rule{.0pt}{2.5ex}
     \\ \hline
        \MAbox{3}{k}
          & 1
          & 1.\overbar{6}
          & [2.\overbar{6},3.\overbar{3}]
          & [4, 4.\overbar{6}]
          & [5.\overbar{6}, 7]
          & [7.\overbar{6}, 9]
          \rule{.0pt}{2.5ex}
   \end{array}
   \end{equation*}

 \subsubsection{Simplices and arbitrary star-regular mesh cells}

 If the mesh cell is not necessarily even, then for
 a regular polygon the optimal interior
 weight is still bounded above by $\MAstar{D}{k}$.
 Based on equations \eqref{WAstar} and \eqref{WAtriplusRecurrence},
 as well as Remark \ref{zhangShuTriangleWeightRemark}
 and Theorem \ref{cubicSimplexWeights},
 we tabulate exact values or bounding intervals for optimal interior weights
 for the triangle ($\MAtri{2}{k}$) and tetrahedron ($\MAtri{3}{k}$).
 Note that to enforce positivity of the retentional for
 the optimal weight, \emph{for a quadratic representation space
 it is sufficient to enforce positivity at the cell center}
 (see Theorem \ref{quadraticRepresentationSpace}), and,
 in the case of a simplex, \emph{for a cubic representation space
 it is sufficient to also enforce positivity at the center
 of each face} (see Theorem \ref{cubicSimplexPoints}).
 When using these weights, take the boundary average 
 as the arithmetic average of face averages (\S \ref{secAffineInvariants}).

   \begin{equation*}
   \begin{array}{l||c|c|c|c|c|c|c|c|}
          k\ (\Uspace=\PkD):
          & 0
          & 1
          & 2
          & 3
          & 4
          & 5
          & 6
          & 7
          \rule{.0pt}{2.5ex}
     \\ \hline
           n=\lfloor k/2\rfloor:
           & 0 & 0
           & 1 & 1
           & 2 & 2
           & 3 & 3
     \\ \hline
        \hline
        \MAintervalk = \frac{(n+1)(n+2)}{2}:
          & 1
          & 1
          & 3 
          & 3 
          & 6 
          & 6 
          & 10
          & 10
          \rule{.0pt}{2.5ex}
     \\ \hline
        \MAstar{2}{k} = \frac{n+1}{2}(\floor{\frac{k+1}{2}}\tplus2)\!
          & 1
          & 1.5
          & 3 
          & 5
          & 6
          & 7.5
          & 10
          & 12
          \rule{.0pt}{2.6ex}
     \\ \hline
        \MAstar{3}{k} = \frac{n+1}{3}(\floor{\frac{k+1}{2}}\tplus3)\!
          & 1
          & \!1.25\!
          & 2.\overbar{6}
          & 3.\overbar{3}
          & 5
          & 6
          & 8
          & 9.\overbar{3}
          \rule{.0pt}{2.6ex}
     \\ \hline
        \hline
        \MAtri{2}{k}
          & 1
          & 1
          & 2 
          & 2.\overbar{2}
          & \tbracket{3.4, 6}
          & \tbracket{3.5 ,6}
          & \tbracket{5\tplus\tfrac{1}{7},10}
          & \tbracket{5.25,10}
          \rule{.0pt}{2.6ex}
     \\ \hline
        \MAtri{3}{k}
          & 1
          & 1
          & 1.\overbar{6}
          & \!1.8\overbar{3}\!
          & \tbracket{2.5\overbar{6},5}
          & \tbracket{2.\overbar{6},6}
          & \tbracket{3\tplus\!\tfrac{9}{14}, 8}
          & \tbracket{3.75, 9.\overbar{3}}
          \rule{.0pt}{2.6ex}
   \end{array}
   \end{equation*}
\end{minipage}}
 \caption{Essential results of Section \ref{CalculationOfWeights}.
 A bounding interval is given where the exact value is unknown.}
 \label{tabulatedWeights}
\end{figure}

 For practical use, we summarize our calculation of interior weights
 in Figure \ref{tabulatedWeights}.

\subsection{Additional definitions for regular polytopes}


To calculate quadrature rules that yield
admissible boundary weights for higher-dimensional mesh cells,
we try to reduce to a one-dimensional problem
using the following symmetry framework.

\begin{definition}[notation and conventions]
  \emph{Let $\Cavg = \avgK$ denote the cell average.
  With star-regular polytopes in mind,
  assume that the canonical mesh cell is centered on the origin.
  Define the \defining{radius} $r$ to be $\nhat\dotp\x$,
  which for a regular polytope is the distance from the origin
  to the closest point on the boundary of the mesh cell.
  We will assume unless stated otherwise (and without loss of generality)
  that the radius of $\K$ is 1.
  Let $\Bavg_r$ denote the average over the boundary
  of the mesh cell $r\K$ rescaled to have radius $r$.}
\end{definition}

\begin{proposition}
 [The cell average is a weighted average of boundary
 averages over rescaled cells]
 \label{cellVsBoundaryAvg}
 Let $\K$ be a star-regular polytope in $D$-dimensional space.
 Then
 \begin{gather}
  \label{CavgFormula}
  \Cavg = \frac{\int_0^1 r^{D-1} \Bavg_r \dr}{\int_0^1 r^{D-1}\dr}
        = D \int_0^1 r^{D-1} \Bavg_r \dr.
 \end{gather}
 Furthermore,
 if the polynomial representation space $\Uspace$ is a subset
 of $\PkD$ then $\Bavg_r$ is a polynomial in $r$ of
 degree at most $k$.
\end{proposition}
\begin{myproof}
  Equation \eqref{CavgFormula} equates weighted averages and is
  justified by observing that
  the thickness of the infinitesimal shell $[r,r+\dr]\cdot\K$ is $\dr$
  and its area is proportional to $r^{D-1}$.
  Since any polynomial is a sum of homogeneous
  polynomials, to justify the final statement
  it is enough to observe that the statement is correct
  for homogeneous polynomials; indeed, if $\U$ is homogeneous
  of degree $k'$
  then so is $\Bavg_r$.
\end{myproof}

\subsection{Linear and quadratic representation spaces}\ 

\begin{theorem}[linear representation space]
Let $\K$ be a star-regular polytope
and let $\Uspace=\PPD^1$.
Assume that the origin is the only point whose orbit
under the isometries of $\K$ is a singleton.
Then the optimal interior weight is $\MAopt=1$
and the positivity-preserving time step is
guaranteed simply by enforcing positivity at the boundary nodes.
\end{theorem}
\begin{myproof}
If $\U$ is a linear function then
$\BAhat(\U)=1$, as can be seen by averaging over the orbit
of $\U$ under the group of isometries of the polytope.
\end{myproof}

\begin{theorem}[quadratic representation space]
\label{quadraticRepresentationSpace}
Let $\K$ be a star-regular polytope, and let $\Uspace=\PPD^2$.
Assume that the origin is the only point whose orbit
under the isometries of $\K$ is a singleton.
Then the optimal time step is guaranteed simply by
enforcing positivity at the boundary nodes and
at the cell center.  Enforcing positivity at
the cell center is equivalent to enforcing positivity
of the \retentional $\RA\opt=\MAopt\Cavg - \Bavg$,
where the optimal interior weight is
$\MAopt=\frac{D+2}{D}$.
%
\end{theorem}
\begin{myproof}
As seen from equation \eqref{CavgFormula},
for any homogeneous quadratic $\U$,
$\BAhat\U=\frac{D+2}{D}$.
If $\U$ is a constant-valued function,
then $\BAhat\U=1$.
So the optimizer cannot be constant. 
So by Theorem \ref{nonEmptyZeroSet}
its zero set must be nonempty.
We can require the optimizer
to be invariant under the isometries of the polytope.
So it must be a positive rotationally invariant quadratic
that has a zero at the origin, i.e.\ a homogeneous quadratic.
The support of the optimal interior sum is therefore restricted to the origin.
So the \retentional is proportional to the value at the cell center.
\end{myproof}

\subsection{Cubic representation spaces}\ 

\begin{theorem}[odd-degree representation spaces for boxes]
\label{oddDegreeForBoxes}
For boxes,
$\MAbox{D}{2n+1}=\MAbox{D}{2n}$, and likewise for other
even polytopes and for spheres.
In particular, for a cubic representation
space, $\MAopt = \frac{D+2}{D}$, and positivity of the interior
sum can be enforced simply by enforcing positivity of the
value at the origin.
\end{theorem}
\begin{myproof}
The isometries of an even polytope by definition include
negation $\x\mapsto -\x$.
Therefore an optimizer $\U\opt\in\PPD^{2k+1}$
symmetric under the isometries of the polytope
lacks odd-degree terms, so $\U\opt\in\PPD^{2k}$.
\end{myproof}

\begin{theorem}[cubic representation space for simplices: optimal retentional points]
Assume that $\Uspace=\PPD^3$.
For any simplex, optimal retentional points are located
at the cell center and at the centers of the faces.
\end{theorem}
\begin{myproof}
\label{cubicSimplexPoints}
  By Figure \ref{mutualConfirmationFig},
  the set of optimal retentional points
  must be zeros of an optimizer and
  can be required to be invariant under the isometries of $\K$.
  Since $\Uspace$ is cubic, this means that
  optimal retentional points can exist only at the center of the cell
  and at the center of its faces.
  (Note that an optimizer candidate uniformly zero on the boundary
  could not have a \outcrowding exceeding zero.)
\end{myproof}

\begin{theorem}[cubic representation space: optimal weight for triangle and tetrahedron]
\label{cubicSimplexWeights}
Assume that $\Uspace=\PPD^3$.
For a triangle the optimal interior weight is
$\MAtri{2}{3}=20/9$,
and the optimal retentional points are located
at the cell center and at the centers of the
edges.
For a tetrahedron the optimal interior weight is
$\MAtri{3}{3}=11/6$,
\end{theorem}
\begin{remark}
  \label{zhangShuTriangleWeightRemark}
 \emph{In \cite{article:ZhangXiaShu12trimesh}, Zhang, Xia, and Shu
  construct a boundary-weighted quadrature rule 
  for triangles.  For representation spaces
  $\PPD^{2n}$ and $\PPD^{2n+1}$, their rule has
  interior weight $\MAinterval^{2n}=\MAinterval^{2n+1}=\frac{(n+1)(n+2)}{2}$
  (see equation \eqref{weightsFor1Dintervals}),
  showing for example that $\MAtri{2}{3}\le 3 = \MAinterval^{3}$.}
\end{remark}
\begin{proofsketch}
  By Theorem \ref{cubicSimplexPoints} and
  Figure \ref{mutualConfirmationFig},
  we may demand an optimizer that is
  invariant under the isometries of $\K$
  and that is zero at the center of the mesh cell and at the centers of its faces. 
  Both for a triangle and for a tetrahedron,
  up to multiplication by a nonzero scalar
  there exists a unique such cubic polynomial.
  This polynomial is definite
  (i.e.\ positive after negating if necessary),
  and evaluating $\BAhat$
  for this polynomial
  yields the values
  $\tfrac{20}9$ for a triangle
  and
  $\tfrac{11}6$ for a tetrahedron.
\end{proofsketch}

\begin{figure}[]
\fbox{\begin{minipage}{\linewidth}
\begin{center}
  \large{\textbf{Key results: optimal retentional positivity points}}
\end{center}
\begin{gather*}
 \begin{matrix} 
  \begin{matrix} 
    \begin{matrix} 
     \begin{matrix}
       \mbox{\includegraphics[width=.22\textwidth]{figures/ex2.0}}
     \end{matrix}
     &
     \begin{matrix}
       \mbox{ \includegraphics[width=.22\linewidth]{figures/ex2.9}}
     \end{matrix}
    \end{matrix}
    \\
    \begin{matrix}
      \mbox{\parbox{.47\textwidth}{
       For quadratic and cubic polynomial representation spaces,
       the set of optimal interior points for a box
       is simply the cell center (Theorems
       \ref{quadraticRepresentationSpace}
       and \ref{oddDegreeForBoxes}).}}
    \end{matrix}
    \\
    \rule{.47\textwidth}{1pt}
    \\
    \begin{matrix}
     \begin{matrix}
       \mbox{\includegraphics[width=.22\textwidth]{figures/ex2.7}}
     \end{matrix}
     &
     \begin{matrix}
       \mbox{ \includegraphics[width=.22\linewidth]{figures/ex2.8}}
     \end{matrix}
    \end{matrix}
    \\
    \begin{matrix}
      \mbox{\parbox{.47\textwidth}{
       For tensor product polynomial spaces, the optimal
       interior points for a box are indicated in
       Remark \ref{tensorProductWeights}
       (see \cite{article:ZhangShu10rectmesh}).}}
    \end{matrix}
  \end{matrix} 
  &
  \begin{matrix} 
    \begin{matrix}
     \begin{matrix}
      \mbox{\includegraphics[width=.22\textwidth]{figures/ex2.3}}
     \end{matrix}
     &
     \begin{matrix}
      \mbox{ \includegraphics[width=.22\linewidth]{figures/ex2.5}}
     \end{matrix}
    \end{matrix}
    \\
    \begin{matrix}
      \mbox{\parbox{.47\textwidth}{
       For quadratic polynomials in a simplex,
       the set of optimal interior points is
       simply the cell center
       (Theorem \ref{quadraticRepresentationSpace};
       contrast with Figure 1 in \cite{article:ZhangShu11}).
       }}
    \end{matrix}
    \\
    \mbox{\parbox{.47\textwidth}{\rule{.47\textwidth}{1pt}\medskip}}
    \\
    \begin{matrix}
     \begin{matrix}
      \mbox{\includegraphics[width=.22\textwidth]{figures/ex2.4}}
     \end{matrix}
     &
     \begin{matrix}
      \mbox{ \includegraphics[width=.22\linewidth]{figures/ex2.6}}
     \end{matrix}
    \end{matrix}
    \\
    \begin{matrix}
      \mbox{\parbox{.47\textwidth}{
       For cubic polynomials in a simplex,
       the set of optimal ``interior'' points consists
       of the cell center and the centers of the faces
       (Theorem \ref{cubicSimplexPoints}).
       }}
    \end{matrix}
  \end{matrix} 
 \end{matrix}
\end{gather*}
\end{minipage}}
\caption{Optimal interior points for standard mesh cells and low-order
  polynomials.
  \emph{Black positivity points are retentional (i.e.\ ``interior'') points.
  Gray points (when shown) are boundary nodes and depend on one's
  choice of a correct boundary quadrature.
  The black points are determined by mesh cell geometry
  and the representation space and are thus independent of
  the choice of gray points.
  Enforcing positivity at the black points automatically enforces positivity
  of the optimal retentional $\RA\opt:=\MAopt\CA-\BA$
  and guarantees that positivity of the cell average is
  maintained for the same minimum time step that one would
  guarantee by enforcing positivity everywhere in the mesh cell
  (assuming that the same wave speed caps are enforced,
  see Section \ref{sec:desing}).
  Positivity is also enforced at boundary nodes, but for a different
  purpose: to ensure that Riemann problems are solved with physical states.
  For a linear representation space, the set of retentional 
  points is empty, and it is enough to enforce
  positivity at the boundary nodes.}}
\label{optimalPositivityPoints}
\hrule
\end{figure}

\medskip
Results for quadratic and cubic polynomials
are summarized in Figure \ref{optimalPositivityPoints}.

\subsection{High-order representation spaces}\ 

\begin{theorem}[weights for 1D interval]
  \label{1Dweights}
  \begin{gather}
    \label{weightsFor1Dintervals}
    \MAinterval^{2n} = \MAinterval^{2n+1}
    = \WA\inv := \frac{(n+1)(n+2)}{2}.
  \end{gather}
  Note that the optimal boundary-weighted quadrature 
  is the Gauss-Lobatto quadrature with $n$ interior points
  and end weight $\WA$.
\end{theorem}
\begin{myproof}
For the canonical mesh cell $\K=\unitinterval$,
the Gauss-Lobatto quadrature
\begin{gather}
  \label{GLQ}
  2\Qavg_{\unitinterval} \U = \WA(\U(0)+\U(1)) + \sum_{i=1}^n \what_i \U(x_i)
\end{gather}
with $n$ interior points
is exact for polynomials of degree at most $2n+1$.
The function 
$ 
  \U\opt(x) = \prod_{i=1}^n (x-x_i)^2
$ 
is zero at all interior points and is of degree $2n$. 
So for the representation spaces $\PPD^{2n}$ and $\PPD^{2n+1}$,
the quadrature rule is exact, $\U\opt$ is a positive
function in the representation
space, and $\U\opt$ is zero on the interior points.
Therefore the hypotheses of Theorem \ref{mutualConfirmation}
are satisfied and the conclusion follows.
\end{myproof}
\begin{remark}
\emph{For 1D mesh cells, Zhang and Shu enforce positivity at the
Gauss-Lobatto quadrature points, which of course implies positivity
of the \retentional \cite{article:ZhangShu10rectmesh}.}
\end{remark}

The previous theorem is a special case of the following theorem.
 \begin{theorem}[optimal interior weight for spherical mesh cells]
 \label{WsphThm}
 \begin{equation}
    \label{WAsph}
    \MAsph{D}{2n}
    = \MAsph{D}{2n+1}
    = \frac{(\lfloor\frac{n}{2}\rfloor+1)(2\lfloor\frac{n+1}{2}\rfloor+D)}{D}.
 \end{equation}
  Furthermore, $\MAboxDk\le\MAsphDk$, with equality precisely when
  $D=1$ or $k<4$.
 \end{theorem}
 \begin{remark}
  \label{tensorProductWeights}
  \emph{For boxes one can also
  construct a boundary-weighted quadrature rule
  by taking the tensor product of a Gauss-Lobatto
  quadrature rule with a quadrature rule used to
  integrate over a face and averaging over all $D$
  such quadrature rules, as is done in
  \cite{article:ZhangShu10rectmesh}
  for the case $D=2$.  The resulting interior weight
  is $\MAopt=\MAintervalk$ independent of $D$.
  This weight and set of quadrature points is optimal
  for a tensor product polynomial space; indeed,
  invoking Corollary \ref{mutualConfirmation},
  a confirming optimizer is 
  $\U\opt(\x):=\prod_{i=1}^D \UoptB{\unitinterval}{k}(x_i)$,
  where $\UoptB{\unitinterval}{k}$ is a unit-interval optimizer.}
  %
  \emph{But for $\PkD$, in case $D>1$ and $k>1$,
  the boundary weight of this theorem
  is an improvement over $\MAintervalk$.}
 \end{remark}

 \begin{proofsketch}
  Recall from \eqref{CavgFormula} that
  $
    \Cavg = D \int_0^1 r^{D-1} \Bavg_r \dr,
  $
  and observe that, for boxes and spheres, $\Bavg_r$ is an even
  polynomial in $r$ of degree less than or equal
  to the degree $k$ of the representation space.
  So we can write
  $\Bavg_r = B(r^2)$,
  where $\mydeg B = \lfloor k/2 \rfloor$.
  We seek a quadrature rule that maximizes
  the weight on the point $r=1$.
  Making the substitution
  $
    t = 1-2 r^2
  $
  shows that $\Cavg$ is given by an integral with a
  Jacobi weight function.

  In case $k=4m$ (or $k=4m+1$) we use a Gauss-Radau quadrature rule
  for Jacobi weight (GRJ) with $m$ interior points, which is exact
  for polynomials of degree at most $2m$
  (see \cite{gautschi00GRJ}).
  In the case of a spherical mesh cell,
  an optimizer that confirms that the boundary weight is maximal is
  $
    \U\opt(r^2) = \prod_{i=1}^m (r^2-r_i^2)^2,
  $
  where the $r_i$ are the interior quadrature points corresponding
  to the interior quadrature points $t_i$ of the GRJ quadrature rule.

  In case $k=4m+2$ (or $k=4m+3$), we use a Gauss-Lobatto quadrature rule
  for Jacobi weight (GLJ) with $m$ interior points, which is exact
  for polynomials of degree at most $2m+1$ (see \cite{gautschi00GL}).
  In the case of a spherical mesh cell,
  an optimizer that confirms that the boundary weight is maximal is
  $
    \U\opt(r^2) = r^2\prod_{i=1}^m (r^2-r_i^2)^2,
  $
  where the $r_i$ are the interior quadrature points corresponding
  to the interior quadrature points $t_i$ of the GLJ quadrature rule.
 \end{proofsketch}

 \begin{theorem}[admissible interior weight for any star-regular mesh cell]
  \label{WAstarThm}
   \begin{equation}
     \label{WAstar}
     \MAtri{D}{k} \le
     \MAstar{D}{k}
     := \frac{(\lfloor \frac{k}{2}\rfloor+1)(\lfloor\frac{k+1}{2}\rfloor+D)}{D}.
   \end{equation}
 \end{theorem}
 \begin{remark}
   \emph{Equality holds only for the case $\MAtri{D}{0} = 1$.
   Recall from Remark \ref{zhangShuTriangleWeightRemark}
   that $\MAtri{2}{2n}\le \MAtri{2}{2n+1} \le \frac{(n+1)(n+2)}{2}$,
   which agrees with estimate \eqref{WAstar} for even-dimensional representation
   spaces and for odd-dimensional representation spaces is slightly better
   than estimate \eqref{WAstar}.}
 \end{remark}
 \begin{proofsketch}
   We want a quadrature rule with such weight on the
   boundary average for
    $
      \Cavg = D \int_0^1 r^{D-1} \Bavg_r \dr.
    $
   Unlike for boxes and spheres, we cannot assume
   that $\Bavg_r$ is an even polynomial.
   For $\PPD^{2n}$ we use a GRJ quadrature
   (see \cite{gautschi00GRJ})
   with $n$ points to a get a quadrature rule for $\Cavg$,
   and for $\PPD^{2n+1}$ we use a GLJ quadrature
   (see \cite{gautschi00GL})
   with $n$ points to get a quadrature rule for $\Cavg$;
   in each case, the resulting end weight is the reciprocal of the
   interior weight appearing in \eqref{WAstar}.
 \end{proofsketch}

 \subsection{Lower bounds for interior weights}

 In the previous sections we computed exact
 values and upper bounds for optimal interior weights.
 In the case of higher-order representation spaces
 where we have merely computed upper bounds,
 we would also like to have a lower bound
 that suggests how much the estimate might be improved.
 In particular, we will obtain bounding intervals that
 prove that
 \begin{align}
 \MAsphDk &= \Omega(k^{2}),
 &\MAboxDk &= \Omega(k^{2}),
   &&\hbox{and} &
 \MAtriDk &= \Omega(k^{2}).
 \label{quadraticBoundaryWeightScaling}
 \end{align}

 Recall that for 1D intervals the optimal weight
 is given by equation \eqref{weightsFor1Dintervals}.
 For boxes and for simplices, we now bootstrap from this result
 to obtain lower bounds for the interior weight.

 \newcommand{\WAboxplus}[1][D]{\WAbox{#1}{k,+}}
 \newcommand{\MAboxplus}[1][D]{\MAbox{#1}{k,-}}
 \subsubsection{Boxes}
 First consider the case of the box $\unitinterval^D$.
 It has $2D$ faces, each of which is a box of dimension $D-1$.
 Each face has $2(D-1)$ sub-faces, each of which is a
 box of dimension $D-2$.
 Suppose that $\U$ is defined on the $x=1$
 face of the box and has interior weight $\MAboxplus[D-1]$
 when restricted to this face.
 Then $\U$ can also be thought of as a function
 on $\unitinterval^D$ and has boundary weight given by
 the recurrence relation
 \begin{gather}
   \label{WAboxplusRecurrence}
   \MAboxplus[D]
   = \frac{1+(D-1)\MAboxplus[D-1]}{D},
 \end{gather}
 as will be shown in detail in Part II of this work \cite{article:Jo12};
 we can take $\MAboxplus[]:=\MAintervalk$.
 By induction on $D$, the recurrence relation \eqref{WAboxplusRecurrence}
 implies that $\MAboxDk = \bigoh(k^2)$,
 so by Theorem \ref{WsphThm}, $\MAboxDk = \Omega(k^2)$.

 \newcommand{\WAtriplus}[1][D]{\WAtri{#1}{k,+}}
 \newcommand{\MAtriplus}[1][D]{\MAtri{#1}{k,-}}
 \subsubsection{Simplices}

 Let $K$ be a regular simplex of dimension $D$.
 Suppose that $\U$ is a polynomial of degree at most $k$
 defined on the face $F$ of $K$ not containing vertex $V$,
 and suppose that $\U$ is positive when restricted to $F$.
 Then $\U$ extends uniquely to a polynomial
 that is homogeneous of degree $k$ when expanded about $V$.
 This yields a recurrence relation for a lower bound
 $\MAtriplus[D]\le\MAtri{D}{k}$,
 \begin{gather}
   \label{WAtriplusRecurrence}
   \MAtriplus[D]
     = \left(\frac{D+k}{D}\right)
       \left(\frac{1+D\frac{D-1}{D+k-1}\MAtriplus[D-1]}{1+D}\right),
 \end{gather}
 as will be shown in detail in Part II of this work \cite{article:Jo12}.
 By induction on $D$, recurrence relation
 \eqref{WAtriplusRecurrence}
 implies that $\MAtriDk = \bigoh(k^2)$.
 Thus, by Theorem \ref{WAstarThm},
 $\MAtriDk = \Omega(k^2)$.


\section{Practical application of positivity limiting}
\label{practicalApplication}
We briefly consider the most prominent issues that arise
when enforcing positivity.
Detailed treatment of practical issues is deferred to
Part II of this work \cite{article:Jo12}.
\begin{figure}[]
\fbox{\begin{minipage}{\linewidth}
\begin{center}
  \large{\textbf{Prescriptions for positivity limiting}}
\end{center}

    \medskip
    \textbf{\emph{Check positivity efficiently.}}
    One can use interval arithmetic to inexpensively confirm cell
    positivity in the vast majority of mesh cells.  
    For nodal DG, choose positivity points to be nodal
    points to avoid additional computational expense.
    See \S \ref{EfficientImplementationOfPositivityChecks}.
    %

    \medskip
    \textbf{\emph{Pad inequalities.}}
    In practice, error in machine arithmetic makes it problematic
    to enforce exact inequalities such as $\rho\ge0$ or $p\ge0$.
    Instead, enforce $\rho\ge\rhoeps$ or $p\ge\peps$ for some
    small $\rhoeps>0$ and $\peps>0$.
    See \cite{WangZhangShuNing12}.


    \medskip
    \textbf{\emph{Use positivity points to improve stability.}}
    Enforcing positivity at additional interior points may improve stability
    (likely depending on the choice of oscillation-suppressing limiters).
    Since values used in volume and boundary quadratures must
    be calculated (or available) anyway, positivity of these values
    can be efficiently enforced, and enforcing positivity at the points
    guarantees that state values used in the numerical method are always positive.
    For modal DG, positivity at any finite set of points
    can be efficiently checked, and it makes sense to include
    the optimal retentional points if they are known.

    \medskip
  \textbf{\emph{Use wave speed desingularization for systems,
    especially for shallow water.}}
    When enforcing positivity in shallow water,
    fluxes need to be calculated with remapped states in order
    to desingularize wave speeds.
    See remarks in Algorithm \ref{alg:scalar} and Algorithm
    \ref{alg:system} and Section \ref{sec:desing}.
    For gas dynamics, wave speed desingularization
    is needed when enforcing positivity of the density
    (although more often it is positivity of the pressure
    that needs to be enforced).

  \medskip
  \textbf{\emph{Estimate an optimal time step after enforcing
    the cell positivity condition.}}
    The factor by which the time step must be shortened determines the
    expense of positivity limiting
    (\S \ref{PositivityPreservingVersusStableTimeStep}).
    Enforcing the cell positivity condition guarantees a minimum time
    step $\Dtpos$ that preserves positivity of the cell average;
    see equation \eqref{weightFunctionalCondB}.
    One can directly calculate the maximum one-stage time step
    $\Dtzero$ that maintains positivity of the cell average
    (\S \ref{StatePositivityIndicators})
    in addition to calculating the maximum stable time step $\Dtstab$
    to determine an optimal safe time step.
    For multistage
    and local
    time stepping, one can use
    $\Dtstab$, $\Dtzero$, and optionally $\Dtpos$
    to maintain and iteratively adjust an estimate of a safe time step
    that is both stable and positivity-preserving.
\end{minipage}}
\caption{Key points of Section \ref{practicalApplication}.}
\label{KeyIssuesInPositivityLimiting}
\end{figure}
We summarize in
Figure \ref{KeyIssuesInPositivityLimiting}.


\def\smin{_{\min}}
\def\smax{_{\max}}
\subsection{Efficient implementation of positivity checks}
\label{EfficientImplementationOfPositivityChecks}

For almost all problems, in the large majority of cells one can
confirm that the solution is positive and will remain so for any
stable time step by using interval arithmetic. Thus, one can
enforce positivity with little additional computational expense
per time step.
For example, in the case of gas dynamics, one can quickly confirm
positivity of the pressure in the vast majority of cells by
checking that $P\smin := \rho\smin\energy\smin-\mom\smax^2/2$
is positive, where $\rho\smin$ and $\energy\smin$ are lower
bounds on density and energy in the cell and $\mom\smax$ is an
upper bound on momentum.

The question is how to obtain a tuple of intervals $[\vu\smin,\vu\smax]$
for the components of the conserved variables. In modal DG one can obtain lower and
upper bounds that hold globally in the mesh cell by using lower
and upper bounds on the polynomial basis functions $\varphi_j$.
This works well, because for smooth solutions the coefficients of
higher-order modes decay rapidly as the mesh is refined. Nodal DG
represents the solution using values at nodes; therefore, one can
easily obtain an interval $[\vu\smin,\vu\smax]$ that applies to
all nodal values.  Thus,
even for nodal DG, checking positivity at all nodes is usually
no more expensive than checking positivity at a point.
%
By Theorem \ref{existenceOfFiniteM},
positivity at the nodes of a nodal DG scheme implies positivity of
a retentional with weight $\WA\le\WAopt$, because
the cell average is a positively weighted sum of values at the nodes
(including boundary nodes)
and so is a boundary-weighted quadrature.

%

\def\dUpara{\dU_\parallel}
\def\dUperp{\dU_\perp}
\def\cTheta{\overbar{\Theta}}
\def\Uavghat{\widehat{\Uavg}}
\subsection{Wave speed desingularization is needed when limiting density}
\label{sec:desing}
  Let $\rho$ denote depth in the shallow water case or
  density in the gas dynamics case.
  After enforcing positivity of $\rho$ at boundary nodes,
  wave speeds need to be de-singularized.
  This can be done by modifying states used to compute fluxes
  when the wave speed exceeds an estimate of the maximum
  wave speed in the physical solution.

  In both shallow water and gas dynamics, desingularization
  of momentum (i.e.\ fluid velocity) is necessary.
  Let $\ucap$ be a fluid speed cap that, for a sufficiently
  refined mesh, is guaranteed to exceed the fluid speed that
  arises in the exact solution.  Then
  an accuracy-respecting
  desingularization map of the momentum (or fluid velocity) is
  the rescaling
  $
    \bu \leftarrow \bu\cdot\Ra(\ucap/\uspd),
  $
  where $\Ra(x)=1$ if $x\ge 1$ and equals $x$ (or the
  spline $x(2-x)$) if $0\le x\le 1$
  (see Section 4.4 of \cite{BolNoeLukMed11}).

  %

  For shallow water, enforcing positivity means enforcing positivity
  of $\rho$ and therefore essentially always entails wave speed desingularization.
  %
  For gas dynamics, in addition to desingularizing
  the fluid speed $\bu$, it may also be necessary to
  desingularize the sound speed $c=\sqrt{\gamma\theta}$,
  e.g.\ by capping the pseudo-temperature $\theta:=p/\rho$.
  %


  More careful study of wave speed desingularization
  is left to future work.

\newcommand{\WAS}[2]{\ensuremath{{\WA_{#2}^{#1,\mathrm{stable}}}}}
\newcommand{\MAS}[2]{\ensuremath{{\MA_{#2}^{#1,\mathrm{stable}}}}}
\newcommand{\WAH}[2]{\ensuremath{{\WA_{#2}^{#1,\mathrm{hom}}}}}
\newcommand{\MAH}[2]{\ensuremath{{\MA_{#2}^{#1,\mathrm{hom}}}}}
\subsection{Cost of positivity limiting}
\label{PositivityPreservingVersusStableTimeStep}


For most problems, for the vast majority of cell updates,
a quick check (\S \ref{EfficientImplementationOfPositivityChecks})
will confirm that positivity limiting is unnecessary.
If local time stepping is used, then the total cost of
positivity limiting will be marginal. But if global time stepping
is used, then the expense of positivity limiting is measured by
the factor by which the time step must be shortened to maintain positivity.
This suggests (1) a careful consideration of how to obtain
tight and reliable wave speed bounds for use in wave speed
desingularization and (2) a comparison of stable and
positivity-preserving time steps.


For (2), in the case of one-dimensional mesh cells
with $k$th order polynomial space, values are known:
an approximate value for the maximum stable time step is given by
$\Dtstab\inv \approx (k+1/2)\cs/\Dx$ if an SSP-RK time-stepping
method of order $k$ is used (see Table 2.2 in \cite{CoShu01});
in comparison, capping \outcrowding by $\MAopt=(n+1)(n+2)/2$ gives
$\Dtpos\inv = \csA\MAopt = (n+1)(n+2)\cs/\Dx$,
where $k$ equals $2n$ or $2n+1$.

\subsection{Isoparametric mesh cells}

For a typical isoparametric mesh cell, in canonical coordinates
the mesh cell is a box or simplex and the representation space is
a polynomial space (e.g.\ $\PkD$). Each
physical mesh cell is the image of such a canonical mesh cell
under a diffeomorphism $\phi$, and the physical representation
space is the push-forward under $\phi$ of the canonical
representation space.
%
%
The outflow positivity limiting framework is defined
in canonical coordinates and can be applied
without modification for isoparametric mesh cells: under
reasonable regularity assumptions on $\phi$, direct application
of positivity limiters in canonical coordinates respects accuracy
and guarantees a minimum positivity-preserving time step
\cite{article:Jo12}.

\section{Conclusion}
\label{Conclusion}

This work consists of two main results: a framework and algorithm
for positivity limiting, and calculations of optimal weights
and positivity points needed by this framework.

\subsection{Outflow positivity limiting framework}
We have developed a framework for preserving positivity
of each cell average based on limiting the rate and amount of material
that can flow out of the cell. This is achieved by 
linear damping of high-order corrections,
remapping boundary node states to limit wave speeds,
and limiting the time step.
In each cell, the high-order corrections are linearly damped
just enough to enforce a physically justified cap $\MA>1$ on the
\mention{boundary crowding} $\BAhat\Alpha(\U):=\frac{\BA\Alpha(\U)}{\CA\Alpha(\U)}$
for all affine functionals $\Alpha$ in a set of state positivity functionals $\dP$,
where $\CA$ is the cell average and $\BA$
is the arithmetic average over all faces
of the average over each face (see Section \ref{affineInvariantDefs});
this condition is satisfied precisely when the \mention{unital retentional}
\begin{gather*}
  \Rhat(\U) = \frac{\MA\CA\U - \BA\U}{\MA-1}
\end{gather*}
satisfies positivity.
One way to enforce positivity of the unital retentional is to
enforce that the solution is positive at appropriately chosen
quadrature points in the mesh cell, as Zhang and Shu have done
\cite{article:ZhShu10,article:ZhangShu11}.
Directly enforcing positivity of the unital retentional is computationally
inexpensive and makes it unnecessary to determine and enforce
positivity at these points.
%
For concreteness, we display positivity-preserving algorithms
for scalar conservation laws and shallow water
(Algorithm \ref{alg:scalar}) and for gas dynamics 
(Algorithm \ref{alg:system}).

\subsection{Optimal weights and positivity points}
Just as point-wise positivity limiting prompts the search for optimal
positivity points, retentional positivity limiting prompts the
need for a reasonably close upper bound $\MA$ on the optimal
interior weight $\MAopt[\K](\Uspace)$.
In Section  \ref{sec:accuracy}, we developed a general framework to estimate 
upper and lower bounds on $\MAopt[\K](\Uspace)$.
In Section \ref{CalculationOfWeights},
we applied this framework to tabulate interior weights for boxes and simplices
for polynomial representation spaces.

In practice, the canonical mesh cell is almost always a cube or
a simplex, and, for problems that are likely to entail positivity
limiting, the representation space is likely to be a space
of polynomials of at most cubic degree.
Assuming these conditions,
the set of optimal interior positivity points is very small, and
the essential take-home result of this work is that
to guarantee the same positivity-preserving
time step as if positivity were enforced everywhere in the mesh cell,
in addition to enforcing positivity (and, for shallow water,
desingularizing wave speeds) at boundary nodes,
\emph{it is sufficient to enforce positivity at the cell center},
except that for a simplex with a cubic representation space
one must also enforce positivity at the centers of the faces.
See Figure \ref{optimalPositivityPoints}.

\bigskip

%

 \begin{figure}[]
   {\smaller
   \begin{gather*}
     \begin{aligned}
       &\calP\subset\realsN && \textrm{set of positive states} \\
       &\dP\subset\setof{\Alpha:\realsN\to\reals}
          &&\textrm{affine posit. functionals} \\
       &K \subset\reals^D&&\textrm{canonical mesh cell} \\
       &\PPD \subset\setof{f:\reals^D\to\reals}&&\textrm{polynomials in $D$ vars.} \\
       &\PkD \subset\PPD
         &&\textrm{polynomials of deg. $\le k$} \\
       &\Uspace\supseteq \PkD&&\textrm{representation space} \\
       &\vol := \volK :=\tight{\intK} 1 &&\textrm{measure of cell} \\
       &\area:= \areaK := \tight{\QintdK} 1 &&\textrm{measure of boundary} \\
       &\C(\U) := \tight{\intK}\U &&\textrm{cell integral} \\
       &\CA(\U) := \Uavg := \tfrac1{\vol}\tight{\intK}\U &&\textrm{cell average} \\
       &\Q = \dQK\subset\dK &&\textrm{boundary quadrature pts.} \\
       &\B(\U)  :=\tight{\QintdK} \Um &&\textrm{boundary integral quad.} \\
       &\BA(\U)  :=\area\inv\tight{\QintdK} \Um &&\textrm{boundary avg. quad.}
     \end{aligned}
     \ \ \ 
     \begin{aligned}
       &\X \subset\realsD
         &&\textrm{set of positivity points} \\
       &\cs/2 >0 &&\textrm{speed cap on $\Q$} \\
       &\csA :=\tfrac{\area}{\vol}\cs&&\textrm{scale-invariant version} \\
       &\h \le \cs\Um&&\textrm{numerical flux} \\
       &\hA :=\tfrac{\area}{\vol}\h&&\textrm{scale-invariant version} \\
       &\Uavg\next := 
         \CA\U-\Dt\BA\hA
          &&\textrm{updated cell average} \\
       &\WA \ge\Dt\csA  &&\textrm{boundary weight} \\
       &\RWA := \CA - \WA\BA
         &&\textrm{retentional} \\
       &\RMA := \MA\CA - \BA &&\textrm{\retentional (rescaled)}\\
       &\BAhat(\U) := \BA(\U)/\CA(\U) &&\textrm{\outcrowding} \\
       &\Uplus[\K] \subset \Uspace\setminus\setof{0} &&
         \textrm{functions positive in $\K$} \\
       &\MAopt[\K] := \sup \BAhat(\Uplus[\K])
         &&\textrm{optimal interior weight} \\
       &\MA :=\WA\inv\ge\MAopt[\K]
         &&\textrm{admissible weight}
     \end{aligned}
   \end{gather*}
   }
 \vskip-2ex
 \caption{
   Key definitions for cell positivity.
  }
 \label{f:defsSummary}
 \end{figure}

\noindent
{\bf Acknowledgements.}
This work was supported in part by NSF grant DMS--1016202.


\bibliographystyle{amsplain}
\bibliography{bibliography/BigBib}

\begin{algorithm}
\caption{Positivity-preserving Euler time-step for scalar conservation laws
  and shallow water.}
\label{alg:scalar}
\begin{algorithmic}
  \STATE
  \emph{For shallow water,} $\U$, $\h$, \emph{and} $\f$
  \emph{represent the first component of} $\vU$, $\vh$, \emph{and} $\vbf$, \emph{respectively.}
  For the canonical mesh cell $\K$ and representation space $\Uspace$
  determine an admissible weight $\MA$ such that 
  $\RAMA\U:=\MA\Uavg-\BA\U>0$  $\forall \U=\sum_i U^i\varphi_i > 0$ 
  (see  Figure \ref{tabulatedWeights}),
  where $\BA\U = \sum_j U^j \BA(\varphi_j)$ is the boundary average
  (see \S \ref{secAffineInvariants}).
  For efficiency, precompute the values $\BA(\varphi_j)$.
\STATE
\STATE
   1. $\forall\K$
   rescale the deviation $\vdU$ from the cell average $\vUavg$
   by a percentage $\dampingcoef[\Y]$,
   damping just enough
   so that $\U>0$ at a set of
   positivity points $\X$ that includes all quadrature points in $\QintdK$
   and 
   so that 
   $\RMA(\U)=\MA\Uavg-\BA\U>0$:
  \begin{gather*}
    \vU := \vU^n, \quad \vUavg:=\CA\vU, \quad \vdU:=\vU-\vUavg, \quad
    \Umin = \min_{\x\in\X} \, \U,
    \\
    \dampingcoef[\X] = \min\setof{1, \, \frac{\Uavg}{\Uavg - \Umin }}, \quad
   \dampingcoef[\RMA] = \min \left\{1, \, \frac{(\MA-1)\Uavg}{\BA(\dU)} \right\}, \quad
       \dampingcoef[\Y]:=\min\left(\dampingcoef[\X],\dampingcoef[\RAMA]\right), \\
    \vU \leftarrow \vUavg + \dampingcoef[\Y] \, \vdU.
  \end{gather*}
  If $\X$ is sufficiently rich (see Figure \ref{mutualConfirmationFig})
  then $\dampingcoef[\X]\le\dampingcoef[\RAMA]$ and 
  computing $\dampingcoef[\RAMA]$ is unnecessary.
\STATE
\STATE 2. Compute a safe time-step $\Dt$:
  \[
        (\Dtzero)\inv = \max_{\hbox{all } K}
           \paren{\frac{\QintdK \h(\Um\!\!,\Up\!\!,\nhat)}{\intK\U}}, \quad
         \Dt\inv = \max\setof{\paren{\zsafety\Dtzero}\inv, (\Dtstab)\inv},
  \]
\STATE  where $\Dtstab$ is the maximum stable time-step
  and $0 < \zsafety < 1$ is a safety factor.
\STATE
\STATE 3. \emph{(for shallow water, not scalar case.)}
          Modify the boundary node
          states $\vUQm$ if necessary in order to desingularize wave
          speeds at boundary nodes (see \S \ref{sec:desing}).
\STATE
\STATE 4. $\forall K$ update the solution:
\[
  \intK\vU^{n+1} \varphi^{j}  = \intK\vU \varphi^{j}
    + \Delta t \QintK \vbf(t^n,\x,\vU) \dotp\nabla\varphi^{j}
    - \Delta t \QintdK\vh(\vUQm, \vUQp, \nhat) \, \varphi^{j}.
\]
\end{algorithmic}
\end{algorithm}

\begin{algorithm}
\caption{Positivity-preserving Euler time-step for gas dynamics
  (cf.\ \S3.2 of \cite{article:ZhShu10}).}
\label{alg:system}
\begin{algorithmic}
  \STATE
    Choose large enough values $\rhoeps>0$ and $\peps>0$ 
    that vanish with machine epsilon.
  \STATE        
  \STATE 1. $\forall \K$ 
         evaluate the solution $\vU^n =: \vU=\vUavg+\vdU$
         at a set of positivity points
         $\X$ which includes the set of boundary nodes $Q$
         (i.e.\ the quadrature points of $\QintdK$)
         to determine the nodal states
         $\vU_\X\supset\vUQm$.
  \STATE
  \STATE 2. $\forall K$ linearly damp the high-order corrections $\vdU$
        by a percentage $\dampingcoef[\rho]$, damping just enough
        so that the density $\rho$ is positive 
        at the points $\X$ and so that the density
        retentional $\RAMA(\rho)=\MA\rhoavg-\BA\rho$ is positive:
        \begin{align*}
        &\dampingcoef[\X] =
          \min\setof{1, \, \frac{\rhoavg-\rhoeps}{\rhoavg - \rhomin }}, \qquad
        \rhomin = \min_{\x\in\X} \rho,
    \\
        &\dampingcoef[\RMA] = \begin{cases}
      1 & \quad \text{if} \quad \BA(\drho) \le \paren{\MA-1}\rhoavg-\rhoeps, \\
      \frac{(\MA-1)\rhoavg-\rhoeps}{\BA\drho} & \quad \text{if} \quad
        \BA(\drho) > \paren{\MA-1}\rhoavg-\rhoeps,
      \end{cases} \\
       &\vU \leftarrow \vUavg + \dampingcoef[\rho] \,\vdU, \qquad
       \vU_\X \leftarrow \vUavg +\dampingcoef[\rho]\,\vdU_\X,  \qquad
       \dampingcoef[\rho]:=\min(\dampingcoef[\X],\dampingcoef[\RAMA]).
        \end{align*}
  \STATE 3. $\forall K$ linearly damp the solution
        just enough so that the pressure is positive at the points $\X$:
        \begin{align*}
          &\dampingcoef[p] := \max\setof{\theta\in[0,1]:
          \,(\forall \x\in \X)
          \,p\paren{\vUavg + \theta\,\vdU(\x)}\ge \peps}, \\
          & \vU  \leftarrow \vUavg +\dampingcoef[p]\,\vdU, \qquad
            \vUQm \leftarrow \vUavg +\dampingcoef[p]\,\vdUQm.
        \end{align*}
\STATE 4. Modify the boundary node states $\vUQm$ if necessary
          in order to desingularize wave speeds at boundary nodes
          (see \S \ref{sec:desing}).
\STATE
\STATE 5. $\forall K$ and for some admissible interior weight $\MA$
        linearly damp $\vdU$ so that  the retentional
        $\RMA\vU:= \MA\vUavg-\QavgdK\vUQm$ is positive:
        \begin{align*}
          \dampingcoef[\Q] &:= \max\setof{\theta\in[0,1]:\,
              p\paren{(\MA-1)\vUavg - \theta\BA\vdUQm}\ge 0}, \\
        \vU &\leftarrow 
          \vUavg + \dampingcoef[\Q] \vdU, \qquad
        \vUQm \leftarrow \vUavg +\dampingcoef[\Q]\,\vdUQm.
        \end{align*}
\STATE 6. Compute a stable, positivity-preserving time-step $\Dt$:
  \begin{gather*}
    \Dt := \zsafety\max_{\hbox{all }\K}\setof{\Dt\in[0,\zsafety\inv\Dtstab]:\,
        p\paren{\intK\vU-\Dt\QintdK\vh\paren{\vUQm,\vUQp}}\ge 0},
  \end{gather*}
   where $\Dtstab$ is the maximum stable time-step
   and $0 < \zsafety < 1$ is a safety factor.
\STATE
\STATE 7. $\forall K$ update the solution:
\[
  \intK\vU^{n+1} \varphi^{j} =
    \intK\vU \varphi^{j} + \Dt\QintK \vbf(t^n,\x,\vU)
    \dotp\nabla\varphi^{j} - \Dt\QintdK\vh\paren{\vUQm, \vUQp, \nhat} \, \varphi^{j}.
\]
\end{algorithmic}
\end{algorithm}

\end{document}